\documentclass[12pt]{article}
\usepackage[margin=1in]{geometry}
\usepackage{setspace}
\usepackage{graphicx}
\usepackage[usenames,dvipsnames,svgnames,table]{xcolor} 

\usepackage{kotex}
\usepackage{graphicx}
\usepackage{epsfig}
\usepackage{amssymb, amsmath, amsthm, mathtools}
\usepackage{verbatim}
\usepackage{harvard} 
\usepackage[affil-it]{authblk}
\usepackage{mathrsfs}
\usepackage{here}
\usepackage{multirow}
\usepackage{makecell}
\usepackage{tikz}
\usepackage{enumerate}
\usepackage[shortlabels]{enumitem}
\usepackage{yfonts}
\usepackage{comment}
\usepackage{longtable}
\usepackage{physics} 
\usepackage[ruled,vlined]{algorithm2e}
\usepackage{bm}
\usepackage{url}

\newcommand{\bea}{\begin{eqnarray*}}
\newcommand{\eea}{\end{eqnarray*}}
\newcommand{\bean}{\begin{eqnarray}}
\newcommand{\eean}{\end{eqnarray}}

\newtheorem{theorem}{Theorem}[section]
\newtheorem{lemma}[theorem]{Lemma}

\newtheorem{corollary}[theorem]{Corollary}
\newtheorem{example}{Example}[section]

\newcommand{\ben}{\begin{enumerate}}
\newcommand{\een}{\end{enumerate}}

\everymath{\displaystyle}
\newcommand{\E}{\mathbb{E}}

\newcommand{\iid}{\stackrel{i.i.d.}{\sim}}

\newcommand{\lra}{\longrightarrow}

\newcommand{\bfx}{ \mathbf{x}}
\newcommand{\bfy}{ \mathbf{y}}

\newcommand{\bfX}{ \mathbf{X}}

\newcommand{\bbN}{ \mathbb{N}}

\newcommand{\bbR}{ \mathbb{R}}
\newcommand{\bbX}{ \bfx_n}

\newcommand{\calD}{\mathcal{D}}

\newcommand{\calM}{\mathcal{M}}

\newcommand{\calX}{ \mathcal{X}}

\newcommand{\bepsilon}{ \boldsymbol{\epsilon}}

\makeatletter
\g@addto@macro\normalsize{%
  \setlength{\abovedisplayskip}{4pt}
  \setlength{\belowdisplayskip}{4pt}
  \setlength{\abovedisplayshortskip}{2pt}
  \setlength{\belowdisplayshortskip}{2pt}
}
\makeatother
\begin{document}
	
\title{Conditional Dirichlet Processes and Functional Condition Models}

\author[1]{Jaeyong Lee}
\author[2]{Kwangmin Lee}
\author[3]{Jaegui Lee}
\author[4]{Seongil Jo}

\affil[1]{Department of Statistics, Seoul National University}
\affil[2]{Department of Big Data Convergence, Chonnam National University}
\affil[3]{SK Innovation}
\affil[4]{Department of Statistics, Inha University}

\maketitle

\abstract{In this paper, we study the conditional Dirichlet process (cDP) when a functional of a random distribution is specified. Specifically, we apply the cDP to the functional condition model, a nonparametric model in which a finite-dimensional parameter of interest is defined as the solution to a functional equation of the distribution.
We derive both the posterior distribution of the parameter of interest and the posterior distribution of the underlying distribution itself. We establish two general limiting theorems for the posterior: one as the total mass of the Dirichlet process parameter tends to zero, and another as the sample size tends to infinity.
We consider two specific models—the quantile model and the moment model—and propose algorithms for posterior computation, accompanied by illustrative data analysis examples.
As a byproduct, we show that the Jeffreys substitute likelihood emerges as the limit of the marginal posterior in the functional condition model with a cDP prior, thereby providing a theoretical justification that has so far been lacking.}

\maketitle

\section{Introduction}
The Dirichlet process \cite{ferguson1973bayesian} is a nonparametric prior distribution over probability distributions. First introduced by Ferguson in 1973, the Dirichlet process has since become a cornerstone of Bayesian nonparametrics due to its elegant mathematical properties and practical applicability. One of its most notable features is conjugacy: if $F$ follows a Dirichlet process and $X_1, \ldots, X_n \mid F \sim F$, then the posterior distribution of $F$ remains a Dirichlet process. Additionally, the marginal distribution of the $X_i$’s follows a P\'olya urn scheme, which is central to posterior inference via Gibbs sampling in models based on the Dirichlet process. Other key theoretical results, such as the stick-breaking construction or Sethuraman’s representation, provide both mathematical insight and practical utility.
With the advancement of Markov chain Monte Carlo (MCMC) methods, the Dirichlet process has been employed in a wide range of Bayesian models. Among them, the mixture of Dirichlet processes (MDP) stands out as a foundational model, widely used for density estimation and model-based clustering.

In this paper, we aim to employ the Dirichlet process as part of the prior distribution for a functional condition model, which will be formally defined later. In this context, the full Dirichlet process is not required; rather, a conditional Dirichlet process, constrained by a functional of the distribution function, is necessary.
In probability theory, the definition of conditional probability is neither intuitive nor entirely mathematically satisfactory. To illustrate this, consider the case where a random probability measure $F$ follows a Dirichlet process. When conditioning on a functional $g(F) = t \in \mathbb{R}$, the conditional probability of the event that $F$ lies in a certain set is not defined as a single number, but rather as a $\sigma(g(F))$-measurable function, where $\sigma(g(F))$ is the $\sigma$-field generated by $g(F)$.
\citeasnoun{chang1997conditioning} and \citeasnoun{tjur1975constructive} have made efforts to address this issue. In this paper, we follow the approach of \citeasnoun{chang1997conditioning}. \citeasnoun{chang1997conditioning} define conditional distributions using the concept of disintegration, which provides a more intuitive definition of conditional distributions. Disintegration decomposes a given measure into a family of measures on the level sets of a functional, where each component measure corresponds to a constant multiple of a conditional probability distribution. This formulation aligns well with the intuitive understanding of conditional distributions.
In this paper, we derive the conditional distribution of the Dirichlet process based on the concept of disintegration.

We now define the functional condition model and discuss its Bayesian formulation. A simplified version of the proposed functional condition model is given by
\begin{equation} \label{simple_model}
\begin{aligned}
    \theta & \sim \pi \\
    F \mid g(F) = \theta & \sim \mathcal{D}_\alpha(dF \mid g(F) = \theta) \\
    x_1, \ldots, x_n \mid F, \theta & \iid F,
\end{aligned}
\end{equation}
where $\theta$ is the parameter of interest and the primary focus of the inference,  $g(F)$ denotes a functional of the distribution function $F$, and $\mathcal{D}_\alpha(dF \mid g(F) = \theta)$ represents the conditional Dirichlet process given $g(F) = \theta$. The base measure $\alpha$ of the Dirichlet process is a finite and non-null measure.

In the literature, most existing works have considered the case where $g(F)$ corresponds to a moment of $F$, leading to the term, moment condition model. In contrast, in this paper, $g(F)$ is not restricted to a moment; it can be any functional of $F$. Accordingly, we refer to our framework as the functional condition model. 

In the existing Bayesian literature, moment condition models have primarily been studied using empirical likelihood and the exponentially tilted empirical likelihood (ETEL) methods. The empirical likelihood is obtained by maximizing the empirical likelihood ratio over the level set determined by $g(F)$. When the functional corresponds to a moment of $F$, nearly explicit expressions for the empirical likelihood can be derived. See \citeasnoun{owen2001empirical} for details.   \citeasnoun{lazar2003bayesian} proposed using the empirical likelihood within Bayesian models in a manner analogous to the parametric likelihood. \citeasnoun{yang2012bayesian} applied empirical likelihood to quantile regression models and established the asymptotic normality of the resulting posterior distribution. \citeasnoun{porter2015bayesian} extended empirical likelihood methods to spatial models. \citeasnoun{chaudhuri2008generalized} applied the Bayesian empirical likelihood framework to generalized linear models. Furthermore, \citeasnoun{chaudhuri2017hamiltonian} proposed a Hamiltonian Monte Carlo method tailored for Bayesian empirical likelihood. 

\citeasnoun{schennach2005bayesian} derived a likelihood function for use in moment condition models by taking the limit of the marginal distribution in a nonparametric model, and termed this likelihood the exponentially tilted empirical likelihood (ETEL). \citeasnoun{chib2018bayesian} investigated the application of ETEL in Bayesian moment condition models and demonstrated that the Bernstein–von Mises theorem holds even under model misspecification. \citeasnoun{chib2022bayesian} further extended the ETEL framework to conditional moment condition models.

The empirical likelihood and ETEL methods offer the advantage of enabling Bayesian inference for parameters of interest in a manner analogous to that in parametric models, without requiring Bayesian computation for infinite-dimensional parameters. This is achieved by constructing the likelihood for the parameter of interest within a semiparametric model.
However, these methods have certain limitations: empirical likelihood is not a true likelihood in the strict sense, and the likelihood derived from ETEL is constructed using a prior distribution that is only employed in the derivation of ETEL itself. Furthermore, the ETEL model for conditional moment condition models \cite{chib2022bayesian} requires the construction of somewhat artificial conditions to approximate the conditional distributions.

In this paper, we derive the posterior distribution and the marginal likelihood for the parameter of interest by employing the Dirichlet process, a widely used prior in nonparametric Bayesian modeling. Furthermore, we demonstrate that inference for the parameter of interest can be conducted in the same manner as in parametric models, using the derived marginal likelihood.

Recently, \citeasnoun{alam2025dpglm} proposed the DPGLM, a semiparametric model for the generalized linear model, where the conditional distribution of the response is specified as a covariate-dependent, exponentially tilted, normalized completely random measure. The DPGLM offers an elegant posterior sampling method for the semiparametric model. 

In this paper, we study the Bayesian analysis of functional condition models using the conditional Dirichlet process as part of the prior distribution. We first derive key properties of the conditional Dirichlet process using the concept of disintegration. Similar to the standard Dirichlet process, the conditional Dirichlet process serves as a conjugate prior. When used as a prior, the conditional Dirichlet process induces a marginal distribution for the observations $X_i$ that takes the form of a P\'olya urn sequence, multiplied by a function of the $X_i$’s.

We also define a general form of the functional condition model and specify its prior distribution using the conditional Dirichlet process, from which we derive the corresponding posterior distribution. The proposed functional condition model not only encompasses the moment condition models that have been widely studied in the literature, but also includes a broader class of models.

The structure of the paper is as follows. Section 2 introduces the concept of disintegration and derives the posterior distribution when the conditional Dirichlet process is used as a prior. Section 3 defines a general form of the functional condition model and derives its posterior distribution. Two theoretical results for the functional condition model are also presented. First, we derive the limiting form of the posterior distribution as the total mass of the conditional Dirichlet process prior tends to zero and define the corresponding Bayesian bootstrap posterior. Second, we establish the asymptotic normality of the posterior distribution as the sample size goes to infinity.  Section 4 presents posterior computation methods and data analyses for three specific models: the quantile model and the moment model. The supplementary material provides additional results for the regression model, including posterior computation, empirical analysis, and detailed proofs of the main theorems.

\section{The Conditional Dirichlet Process} 
In this section, we begin with a brief introduction to disintegration. For further details, see \citeasnoun{chang1997conditioning}.

Suppose $\mathcal{X}$ and $\mathcal{T}$ are measurable spaces, $\lambda$ and $\mu$ are $\sigma$-finite measures on $\mathcal{X}$ and $\mathcal{T}$, respectively, and $T : \mathcal{X} \to \mathcal{T}$ is a measurable function. The measure $\lambda$ is said to admit a $(T, \mu)$-disintegration $\{ \lambda_t : t \in \mathcal{T} \}$ if the following three properties hold:
\begin{enumerate}[(i)]
\item (concentration property) Each $\lambda_t$ is a $\sigma$-finite measure on $\mathcal{X}$, and
$$
\lambda_t\{T \neq t\} = 0 \quad \text{for } \mu\text{-almost all } t.
$$
\item (measurability) For every nonnegative, real-valued function $f$ on $\mathcal{X}$, the map
$$
t \mapsto \lambda_t f := \int f(x)  \lambda_t(dx)
$$
is measurable.
\item (iterative integration) For every nonnegative, real-valued function $f$ on $\mathcal{X}$,
$$
\int_\mathcal{X} f(x)  \lambda(dx) = \int_\mathcal{T} \left( \int_{[T=t]} f(x) , \lambda_t(dx) \right) \mu(dt).
$$
\end{enumerate}

Note that in this definition, $\lambda$ need not be a probability measure. If $\lambda$ is a probability measure, then ${ \lambda_t }$ defines constant multiples of the regular conditional probability of $\lambda$.

\citeasnoun{chang1997conditioning} provides a set of conditions for the existence of the disintegration (see Theorem 1). 
A particularly convenient set of sufficient conditions for the existence of disintegration is as follows. 
If $\mathcal{X}$ is a Polish space, $\mathcal{T}$ is a Euclidean space equipped with the Lebesgue measure $\mu$, and $\lambda \circ T^{-1} \ll \mu$, then $\lambda$ admits a unique $(T, \mu)$-disintegration. 

Suppose a random distribution $F$ on a set $\mathcal{X}$ follows a Dirichlet process with parameter $\alpha$; that is,
$
F \sim \mathcal{D}_\alpha,
$
where $\alpha$ is a non-null, finite measure on $\mathcal{X}$, and can be written as
$
\alpha = A \bar{\alpha}$, with $A = \alpha(\mathcal{X})$ and $\bar{\alpha} = \alpha / A.$
Let
$
g: \mathcal{M}(\mathcal{X}) \to \mathbb{R}^k
$
be a functional of $F$, where $\mathcal{M} = \mathcal{M}(\mathcal{X})$ denotes the set of all probability measures on $\mathcal{X}$.
To define the conditional Dirichlet process (cDP), we begin by introducing the following assumption on the functional of the distribution. This assumption is necessary for the existence theorem of disintegration that we will use.\\*

\noindent{\bf A1} When $F \sim \mathcal{D}_\alpha$, the functional  $g(F)$ has a density function  \(h(\cdot: g, \alpha)\) with respect to the Lebesgue measure \(\mu\) 
on \(\mathbb{R}^k\) and  \(h(\cdot: g, \alpha)\) is positive for every value in the range of $g$.   \\*

Condition {\bf A1} is satisfied by a broad class of functionals.  
For example, if \( g \) is defined through a moment condition  
\[ g(F) := \int \tilde{g}(x) \, F(dx), \]  
where \( \tilde{g}: \mathcal{X} \rightarrow \mathbb{R} \), and if \( \alpha \circ \tilde{g}^{-1} \) is not a point measure, then \( g(F) \) possesses a density function  \cite{ghosal2017fundamentals} and thus satisfies {\bf A1}.  
Another example is when \( g(F) = \int_D dF \) for some measurable subset \( D \subset \mathcal{X} \);  in this case, \( g(F) \sim \text{Beta}(\alpha(D), \alpha(D^c)) \), and hence \( g(F) \) has a known density.

Under condition {\bf A1}, the Dirichlet process \( \mathcal{D}_\alpha \) admits a unique \( (g, \mu) \)-disintegration \( \{ \mathcal{D}_{\alpha, \xi} : \xi \in \mathbb{R}^k \} \), where each \( \mathcal{D}_{\alpha, \xi} \) is a finite measure on \( \mathcal{M}(\mathcal{X}) \) supported on the fiber
\[
g^{-1}(\xi) := \{ F \in \mathcal{M}(\mathcal{X}) : g(F) = \xi \}, \quad \xi \in \mathbb{R}^k.
\]
The disintegration theorem implies that for any non-negative measurable function \( f : \mathcal{M}(\mathcal{X}) \to \mathbb{R} \), we have
\begin{align} \label{eq:DP-disintegration}
\int_{\mathcal{M}} f(F) \, \mathcal{D}_\alpha(dF) 
&= \int_{\mathbb{R}^k} \int_{g^{-1}(\xi)} f(F) \, \mathcal{D}_{\alpha, \xi}(dF) \, \mu(d\xi) \nonumber \\
&= \int_{\mathbb{R}^k} \int_{g^{-1}(\xi)} f(F) \, \mathcal{D}_\alpha(dF \mid g(F) = \xi) \, h(\xi; g, \alpha) \, \mu(d\xi),
\end{align}
where \( \mathcal{D}_\alpha(dF \mid g(F) = \xi) \) denotes the conditional distribution of \( F \sim \mathcal{D}_\alpha \) given \( g(F) = \xi \).

The above expression can be rewritten more concisely as
\begin{align*}
\mathcal{D}_\alpha(dF) 
&= \mathcal{D}_{\alpha, \xi}(dF) \, \mu(d\xi) \\
&= \mathcal{D}_\alpha(dF \mid g(F) = \xi) \, h(\xi; g, \alpha) \, \mu(d\xi).
\end{align*}
From this, we obtain the following characterizations:
\begin{align*}
\mathcal{D}_\alpha(dF \mid g(F) = \xi) &= \frac{\mathcal{D}_{\alpha, \xi}(dF)}{\mathcal{D}_{\alpha, \xi}(\mathcal{M})}, \\
h(\xi; g, \alpha) &= \mathcal{D}_{\alpha, \xi}(\mathcal{M}) = \mathcal{D}_{\alpha, \xi}(g^{-1}(\xi)).
\end{align*}

In informal notation,  we often write $\calD_\alpha(dF) I(g(F) = \xi)$ to denote the Dirichlet process restricted to the fiber $(g)^{-1}(\xi)$. However, this expression is not mathematically well-defined,  since in most cases of interest,
\[\calD_\alpha( F: g(F) = \xi ) = 0, ~ \forall \xi. \]
In contrast, the measure  $\calD_{\alpha, \xi}$ that appears in the  $(g, \mu)$-disintegration of $\calD_\alpha$ is mathematically well-defined and is supported on  the fiber $(g)^{-1}(\xi)$.

Consider the model  
\begin{align}\label{eq:md1} 
\begin{split}
F & \sim \mathcal{D}_\alpha, \\
x_1, \ldots, x_n \mid F & \iid F.
\end{split} 
\end{align}
Let \( Q \) denote the joint distribution of \( F \) and \( \mathbf{x}_n = (x_1, \ldots, x_n) \). Then \( Q \) can be expressed as
\begin{align*}
Q(d\mathbf{x}_n, dF) 
&= \prod_{i=1}^n F(dx_i) \, \mathcal{D}_\alpha(dF) \\
&= \mathcal{D}_{\alpha + \sum_{i=1}^n \delta_{x_i}}(dF) \, \text{Polya}_\alpha(d\mathbf{x}_n),
\end{align*}
where 
\( \text{Polya}_\alpha(d\mathbf{x}_n) \) denotes the distribution of the Pólya urn sequence.  

Now consider the nonparametric model in which \( F \) is subject to a functional constraint $
g(F) = \xi \in \mathbb{R}^k.$

We place the prior \( \mathcal{D}_\alpha(\cdot \mid g(F) = \xi) \) on \( F \), restricting \( F \) to lie in the fiber \( g^{-1}(\xi) \). The resulting functional condition model is given by
\begin{equation} \label{eq:basic_model} 
\begin{aligned}
F & \sim \mathcal{D}_\alpha(\cdot \mid g(F) = \xi), \\
x_1, \ldots, x_n \mid F & \iid F.
\end{aligned}
\end{equation}

Below, we derive the posterior distribution under model \eqref{eq:basic_model}, for which we first present the following lemma whose proof is  provided in the supplementary material.

\begin{lemma} \label{lm1}  
Suppose condition {\bf A1} holds. Under model \eqref{eq:basic_model}, we have  
\begin{align*}
    \prod_{i=1}^n F(dx_i) \, \mathcal{D}_{\alpha}(dF \mid g(F) = \xi)
    &= \text{Polya}_{\alpha}(d\mathbf{x}_n) \, \frac{h(\xi; g, \alpha + nF_n)}{h(\xi; g, \alpha)} \\
    &\phantom{00} \times \, \mathcal{D}_{\alpha + nF_n}(dF \mid g(F) = \xi), \quad \mu\text{-a.a. } x,
\end{align*}

where \( F_n := \frac{1}{n} \sum_{i=1}^n \delta_{x_i} \), and \( \delta_x \) denotes the point mass  at \( x \).
\end{lemma}

\begin{theorem} \label{th1} 
Suppose condition {\bf A1} holds and $(F, \bfX_n)$ follow model \eqref{eq:basic_model}. Then,
\begin{enumerate}[(a)]
\item[(i)] \( F \mid x_1, \ldots, x_n \sim \mathcal{D}_{\alpha + nF_n}(dF \mid g(F) = \xi) \);
\item[(ii)] \( x_1, \ldots, x_n \sim \dfrac{h(\xi; g, \alpha + nF_n)}{h(\xi; g, \alpha)} \, \text{Polya}_\alpha(d\mathbf{x}_n) \).
\end{enumerate}
\end{theorem}

Theorem \ref{th1} follows directly from Lemma \ref{lm1}. 

\begin{example} 
In this example, we consider the functional condition model \eqref{eq:basic_model} with
\[
g(F) := F(B) = p \in [0, 1],
\]
for some fixed set \( B \subset \mathcal{X} \). That is,
\begin{equation}\label{eq:quantile_model1}
\begin{aligned}
F &\sim \mathcal{D}_{\alpha + nF_n}(dF \mid F(B) = p), \\
x_1, \ldots, x_n \mid F &\iid  F.
\end{aligned}
\end{equation}

In this case, the distribution of \( F(B) \) is known when \( F \sim \mathcal{D}_\alpha \).  
By applying Theorem \ref{th1}, we obtain the following:
\begin{equation}\label{eq:quantile_post1}
\begin{aligned}
F \mid x_1, \ldots, x_n &\sim \mathcal{D}_{\alpha + nF_n}(dF \mid F(B) = p), \\
(x_1, \ldots, x_n) &\sim \frac{\mathrm{Beta}(p; \alpha(B) + nF_n(B), \alpha(B^c) + nF_n(B^c))}{\mathrm{Beta}(p; \alpha(B), \alpha(B^c))} \, \mathrm{Polya}_n(d\mathbf{x}_n) \\
&= \frac{[\alpha(\mathcal{X})]_n}{[\alpha(B)]_{nF_n(B)} \, [\alpha(B^c)]_{nF_n(B^c)}} \, p^{nF_n(B)} (1 - p)^{nF_n(B^c)} \, \mathrm{Polya}_n(d\mathbf{x}_n),
\end{aligned}
\end{equation}
where \( \mathrm{Beta}(x; a, b) \) denotes the density of the \( \mathrm{Beta}(a, b) \) distribution evaluated at \( x \), and \( [a]_k := a(a+1)\cdots(a+k-1) \) denotes the rising factorial.

In particular, if we assume \( p = \tfrac{1}{2} \), then
\begin{equation} \label{eq:quantile_marginal1}
(x_1, \ldots, x_n) \sim \frac{[\alpha(\mathcal{X})]_n}{[\alpha(B)]_{nF_n(B)} \, [\alpha(B^c)]_{nF_n(B^c)}}  \, \mathrm{Polya}_n(d\mathbf{x}_n).
\end{equation}
When \( n = 1 \), by the tail-free property of the Dirichlet process, we know
\[
x_1 \sim \tfrac{1}{2} \, \overline{\alpha|_B}(dx) + \tfrac{1}{2} \, \overline{\alpha|_{B^c}}(dx),
\]
where \( \alpha|_B \) denotes the restriction of the measure \( \alpha \) to the set \( B \), and \( \overline{\alpha|_B} \) is its normalized probability measure. 
We recover the same result by expanding equation \eqref{eq:quantile_marginal1}.

For general \( n \in \mathbb{N}\), we can obtain a similar expression for the marginal distribution of \( (x_1, \ldots, x_n) \):
\begin{align*}
&\frac{[\alpha(\mathcal{X})]_n}{[\alpha(B)]_{nF_n(B)} \, [\alpha(B^c)]_{nF_n(B^c)}} \cdot \left(\tfrac{1}{2}\right)^n \, \mathrm{Polya}_n(d\mathbf{x}_n) \\
&= \left(\tfrac{1}{2}\right)^n \sum_{i_1 = 0}^1 \cdots \sum_{i_n = 0}^1 \overline{\alpha|_{B^{i_1}}}(dx_1) \, \overline{(\alpha + \delta_{x_1})|_{B^{i_2}}}(dx_2) \cdots \overline{(\alpha + \sum_{j=1}^{n-1} \delta_{x_j})|_{B^{i_n}}}(dx_n),
\end{align*}
where \( B^1 = B \) and \( B^0 = B^c \).
\end{example}

\section{The Functional Condition Model} 

\subsection{The model} 
In this section, we propose a general formulation of the functional condition model and derive its posterior distribution. 
Suppose we observe $n \in \bbN$ pairs of the responses and covariates,  $(x_i, y_i)$, for $i=1,2, \ldots, n$, where the $x_i$'s are nonrandom covariates and the $y_i$'s are random variables. 
The distributions of the pairs $(x_i, y_i)$
are defined through the equation 
\[ \epsilon_i = t(y_i, x_i, \xi) \iid F, ~ i=1,2, \ldots, n, \]
where  $t$ is an $\bbR^{d_\epsilon}$-valued function,   $\xi \in \bbR^k$ is the model parameter, $F$ is an unknown distribution on
$\calX \subset \bbR^{d_\epsilon}$ and $d_\epsilon, k \in \bbN$. 

We assume that the parameter $\xi$ satisfies a functional condition involving the distribution $F$ 
\begin{equation}\label{eq:1}
g(F, \xi)=0 ,
\end{equation}
where $g: \calM \times \bbR^k \to \bbR^{k_1}$ with $k_1 \in \bbN$. 
In some cases,  the parameter $\xi$ is partitioned as  $\xi = (\theta, \nu)$ and the model imposes an additional constraint $\nu = 0$.  
In summary, the functional condition model we propose is as follows. 

\noindent {(Functional condition model)}
\begin{equation}\label{eq:model}
\begin{aligned}
\epsilon_i  = t(y_i, x_i, \xi) & \iid F, ~ i=1,2, \ldots, n, \\
g(F, \xi) & = 0, \\
\nu & = 0. 
\end{aligned} 
\end{equation}
The functional condition model \eqref{eq:model} encompasses a broad class of models.
Depending on the specific application, certain components of model \eqref{eq:model} may be omitted.
The following basic model covers the majority of cases and is relatively easier to work with.

The functional condition model \eqref{eq:model} includes a large class of models. 
Depending on models, parts of the model \eqref{eq:model} can be omitted. The following basic model covers majority of cases and is relatively easier to work with. \\
\noindent (Basic functional condition model) 
\begin{equation}\label{eq:standard_model}
\begin{aligned}
y_1, \ldots, y_n  & \iid F, \\
g^*(F) & = \theta,
\end{aligned} 
\end{equation}
where $g^*: \mathcal{M} \to \mathbb{R}^k.$

The model \eqref{eq:standard_model} does not include nonrandom covariates $x_i$ or the fixed parameter $\nu$. The observations $y_i$ are assumed to follow the distribution $F$, and the functional condition directly determines the parameter of interest, $\theta$.

\subsection{The prior and posterior} 
In this subsection,  we derive the posterior distribution of model \eqref{eq:model}.   
Let $\bfx_n = (x_1, \ldots, x_n)$ and $\bfy_n = (y_1, \ldots, y_n)$.  For each $\bfx_n$, $\bfy_n$ and $\xi$,  define
\[ \bepsilon_n = (\epsilon_1, \ldots, \epsilon_n) = t(\bfy_n, \bfx_n, \xi) = (t(y_1, x_1,  \xi), \ldots,  t(y_n, x_n, \xi) ).\]
Let $\bepsilon^* = (\epsilon_1^*, \ldots, \epsilon^*_{k_\epsilon})$ be the distinct values of $\bepsilon_n$ and let $\Pi(\bepsilon_n)$ denote the 
partition of $[n] =\{1,2, \ldots, n\}$ induced by the equivalence relation
\[ i \sim j \iff \epsilon_i = \epsilon_j, ~ \forall i, j \in [n].\]
To derive the posterior, we impose the following assumption. \\
\noindent{\bf A2}  $\Pi(\bepsilon_n) = \Pi(\bfy_n)$,  for all $\xi$ and $(\bfx_n, \bfy_n)$.

For the prior on $(F, \xi)$, we impose 
\[\pi(d\xi) \propto \pi(\theta)d\theta \delta_0(d\nu), \]
where $\pi(\theta)$ is the prior density of $\theta$. 
We place the conditional Dirichlet process prior for $F$ given the functional condition $g(F, \xi) = 0$,
\[ F \sim \mathcal{D}_\alpha(dF | g(F, \xi) = 0), \]
where $\alpha$ is a finite nonnull measure on $\mathcal{X}$. 
Thus, the prior  of model \eqref{eq:model} is 
\begin{equation} \label{eq:prior} 
\pi(dF, d\xi )  = \calD_\alpha (dF | g(F, \xi) = 0 ) \pi(\theta) d\theta \delta_0(d\nu). 
\end{equation}
In summary, the functional condition model with conditional DP prior is given below.  

\noindent {(Functional condition model with conditional DP prior)}
\begin{equation}\label{eq:model_and_prior}
\begin{aligned}
\epsilon_i  = t(y_i, x_i, \xi) & \iid F, ~ i=1,2, \ldots, n, \\
\xi & \sim \pi(\theta) d\theta \delta_0(d\nu), \\
F & \sim \mathcal{D}_\alpha(dF| g(F, \xi) = 0). 
\end{aligned} 
\end{equation}
To derive the posterior, we assume that $\bar{\alpha}$ has a probability density $a(\cdot)$.  Let $b(y: x,\xi)$ be the density of $y$ 
when $\epsilon \sim a(\cdot)$.  Since $\epsilon = t(y, x, \xi)$, 
\[ b(y: x, \xi) = a(t(y, x, \xi)) | t'(y,  x, \xi)|,   \]
where  $ t'(y, \theta, x) = \frac{\partial}{\partial y}t(y, \theta, x)$.  If the function $t$ does not depend on $\xi$,  the density $b$ of $y$ does not 
depend on $\xi$ either. 

\vspace{-3ex}
\begin{theorem} (posterior) \label{th:post_cDP} 
Given model \eqref{eq:model_and_prior} and under assumptions {\bf A1} and {\bf A2}, we obtain the following results. \\
\noindent (i) (joint posterior) 
\begin{equation}\label{eq:post_cDP}
\begin{aligned} 
\pi(dF, d\theta | \bfy_n ) & \propto \calD_{\alpha + n F_{n, (\theta, 0)} } (dF| g(F, (\theta, 0)) = 0 ) \nonumber\\
&\phantom{00}\times \prod_{i=1}^k  b(y_i^*: (\theta, 0), x_i^*)   \frac{h(0: g(F, (\theta, 0)), \alpha + n F_{n, (\theta, 0)})} { h(0: g(F, (\theta, 0)), \alpha)}  \pi(\theta)d\theta, \\
\end{aligned}
\end{equation}
\noindent (ii) (conditional posterior of $F$) 
\begin{equation} \label{eq:posterior_conditional_F}
\pi(dF| \bfy_n, \theta) =  \calD_{\alpha + n F_{n, (\theta, 0)} } (dF| g(F, (\theta, 0)) = 0 ). 
\end{equation}
\noindent (iii)  (marginal posterior of $\theta$) 
\begin{equation} \label{eq:post_marginal_theta}
\pi(d\theta | \bfy_n )  \propto   \prod_{i=1}^k  b(y_i^*: (\theta, 0), x_i^*)   \frac{h(0: g(F, (\theta, 0)), \alpha + n F_{n, (\theta, 0)})} { h(0: g(F, (\theta, 0)), \alpha)}  \pi(\theta)d\theta. 
\end{equation}
\end{theorem} 

The posterior \eqref{eq:post_cDP} reveals that factor $\prod_{i=1}^k  b(y_i^*: (\theta, 0), x_i^*)$, although originating from the prior, does not vanish as $n \to \infty $.  
We recommend  using a modified posterior in which this factor is removed: 
\begin{equation} \label{eq:post_modified}
\pi^*(dF, d\theta | \bfy_n ) \propto \calD_{\alpha + n F_{n, (\theta, 0)} } (dF| g(F, (\theta, 0)) = 0 )  \frac{h(0: g(F, (\theta, 0)), \alpha + n F_{n, (\theta, 0)})} { h(0: g(F, (\theta, 0)), \alpha)}   \pi(\theta)d\theta.
\end{equation}
The conditional posterior of $F$ of the modified version, given $\theta$, is the same as original posterior \eqref{eq:posterior_conditional_F}, while the marginal posterior of $\theta$ of the modified version is 
\begin{equation}\label{eq:post_modified_marginal_theta}
    \pi^*(d\theta | \bfy_n )  \propto  \frac{h(0: g(F, (\theta, 0)), \alpha + n F_{n, (\theta, 0)})} { h(0: g(F, (\theta, 0)), \alpha)}  \pi(\theta)d\theta. 
\end{equation}

We now apply Theorem \ref{th:post_cDP} to the basic model \eqref{eq:basic_model}. 
The basic model with cDP prior can be summarized as follows. \\
\noindent {(Basic functional condition model with conditional DP prior)}
\begin{equation}\label{eq:basic_model_and_prior}
\begin{aligned}
y_i & \iid F, ~ i=1,2, \ldots, n, \\
\theta & \sim \pi(\theta) d\theta, \\
F & \sim \mathcal{D}_\alpha(dF| g^*(F) = \theta). 
\end{aligned} 
\end{equation}

\begin{corollary} (posterior of the basic model) \label{cor:post_basic_cDP} 
Given model \eqref{eq:basic_model_and_prior} and under assumptions {\bf A1} and {\bf A2}, we obtain the following results. 
\begin{enumerate}[(a)]
\item[(i)] (joint posterior) 
\begin{equation}\label{eq:post_basic_cDP}
\begin{aligned} 
\pi(dF, d\theta | \bfy_n ) & \propto \calD_{\alpha + n F_{n} } (dF| g^*(F) = \theta)   \frac{h(\theta: g^*(F) , \alpha + n F_{n})} { h(\theta: g^*(F), \alpha)}  \pi(\theta)d\theta 
\end{aligned}
\end{equation}
where $F_{n} = \frac{1}{n} \sum_{i=1}^n \delta_{y_i}$.
\item[(ii)] (conditional posterior of $F$) 
\begin{equation} \label{eq:posterior_basic_conditional_F}
\pi(dF| \bfy_n, \theta) =  \calD_{\alpha + n F_{n}  } (dF| g^*(F) = \theta ). 
\end{equation}
\item[(iii)] (marginal posterior of $\theta$) 
\begin{equation} \label{eq:post_basic_marginal_theta}
\pi(d\theta | \bfy_n )  \propto      \frac{h(\theta: g^*(F), \alpha + n F_{n} )} { h(\theta: g^*(F), \alpha)}  \pi(\theta)d\theta. 
\end{equation}
\end{enumerate}
\end{corollary} 


\subsection{The limit of the posterior when $\alpha(\calX) \to 0$}

In this subsection, we consider the limit of the modified posterior  \eqref{eq:post_modified_marginal_theta} for $\theta$  as $A= \alpha(\calX)  \to 0$. We require the following assumptions for the results. 

\begin{enumerate}
    \item[(B1)]
    For all $\theta \in \Theta$, the density $h(\cdot: g(F, (\theta, 0)),  n F_{n, (\theta, 0)} )$ exists. Furthermore, 
as $A \to 0$, we have 
\[ h(0: g(F, (\theta, 0)), \alpha + n F_{n, (\theta, 0)} ) \to h(0: g(F, (\theta, 0)),  n F_{n, (\theta, 0)} ).\]
    \item[(B2)] We can write $h(0: g(F, (\theta, 0)), \alpha  )$ as a product of two factors,
    \[h(0: g(F, (\theta, 0)), \alpha  ) = h_1(\theta,\alpha) h_2(A),\]
    where $h_1(\theta,\alpha)$ is a function that depends on $\theta$ and $\alpha$, and $h_2(A)$ is a function that depends only on $A$. Furthermore, there exists a function $\bar{h}(\theta)$, positive for all $\theta$, such that as $A \to 0$,
    $ h_1(\theta,\alpha) \to \bar{h}(\theta), ~ \forall \theta.$
    
    \item[(B3)] There exists an integrable function $G(\theta)$ such that 
    \bea 
     \Big| \frac{h(0: g(F, (\theta, 0)), \alpha + n F_{n, (\theta, 0)} )}{ h_1(\theta, \alpha  ) }  \pi(\theta)\Big|\le G(\theta).
    \eea 
\end{enumerate}

Condition (B1) is a natural condition stating that, as $A \to 0$, the effect of $\alpha$ disappears from the density of $g(F, (\theta, 0))$. 
Condition (B2) is a technical condition introduced to characterize the limiting behavior of the density \( h(0; g(F, (\theta, 0)), \alpha) \) as \( A \to 0 \), up to a constant. Without accounting for the constant, the denominator may tend to zero as \( A \to 0 \), making it difficult to analyze the limit. By factoring out a function \( h_2(A) \) that depends only on \( A \), we isolate the \( \theta \)-dependent part, \( h_1(\theta, \alpha) \), and ensure that it converges to a nonzero limit \( \bar{h}(\theta) \). 

Condition (B3) is a sufficient condition for ensuring the integrability of the limiting distribution. 
Under these conditions, Theorem~\ref{thm:asymptoticA} provides the limiting distribution of the posterior, and its proof is given in the supplementary material.

\vspace{-2ex}
\begin{theorem}\label{thm:asymptoticA}
Suppose conditions B1 - B3 hold and 
 \[ z_A = \int   \frac{h(0: g(F, (\theta, 0)), \alpha + n F_{n, (\theta, 0)} )}{ h(0: g(F, (\theta, 0)), \alpha  ) } \pi(\theta) d\theta <\infty. \] Then, as $A \to 0$,
\bean\label{eq:convergenceA} 
\pi(\theta\mid \bfy_n) \lra \frac{h(0: g(F, (\theta, 0)),  n F_{n, (\theta, 0)} )}{z_0 \bar{h}(\theta)} \pi(\theta),~\forall \theta,
\eean  
where
 \[z_0 = \int \frac{h(0: g(F, (\theta, 0)),  n F_{n, (\theta, 0)} )}{ \bar{h}(\theta)} \pi(\theta)d\theta<\infty.\] 
\end{theorem}

A similar limiting distribution can be obtained for the posterior  \eqref{eq:post_basic_marginal_theta} of the basic model, where 
$h(0; g(F, (\theta, 0), \alpha + n F_{n, (\theta, 0)}))$
and $h(0; g(F, (\theta, 0), \alpha))$
are replaced by
$h(\theta; g^*(F), \alpha + n F_{n, (\theta, 0)})$
 and $h(\theta; g^*(F), \alpha)$,
respectively. 


\subsection{Limit of the posterior when $n \to \infty$}
In this subsection, we investigate the asymptotic normality of the modified marginal posterior \eqref{eq:post_modified_marginal_theta} of $\theta$. Let $\pi^*(\theta\mid \bfx_n) \propto \exp(f_n(\theta)) \pi(\theta)$, where
$f_n: \Theta \subseteq \bbR^p \to  \bbR$ is defined as
\bea 
f_n(\theta):= \log \frac{h(0:g(F,(\theta,0)), \alpha + n F_{n, (\theta, 0)} )}{h(0: g(F,(\theta,0)),\alpha)}.
\eea 

We impose the following assumptions on $f_n(\theta)$ and $\pi(\theta)$ for the asymptotic normality. 
\begin{enumerate}[(C1)]
    \item[(C1)]
    Let $\theta_n\in \bbR^p$ be a sequence satisfying $\theta_n \to \theta_0\in \Theta$ where $\theta_0$ is  the true parameter value of $\theta$. Suppose that $ f_n(\theta)$ can be expressed as 
    \bea 
    f_n(\theta) = -\frac{n}{2} (\theta-\theta_n)^T H_n (\theta-\theta_n) + C_n + r_n(\theta),
    \eea 
    where $C_n$ is a real sequence, $H_n \in \mathbb{R}^{p \times p}$ is a symmetric matrix such that $H_n \to H_0$ for some  positive definite matrix $H_0 \in \mathbb{R}^{p\times p}$,  and $r_n : \mathbb{R}^p \to \mathbb{R}$ is a remainder term. Furthermore, $r_n$ satisfies the following properties:
    \bean 
    |r_n(x/\sqrt{n}+\theta_n) - r_n(\theta_n) |  &\to& 0, ~\textup {for every fixed }x\in \bbR^p, \label{eq:rn1}\\
    |r_n(x/\sqrt{n}+\theta_n) - r_n(\theta_n) | &\le& \frac{c_0 ||x||^3}{2\sqrt{n}} + 
     x^T c_n+
    d_n,~ x\in B_{\sqrt{n}\epsilon}(0),\label{eq:rn2}
    \eean 
    where $c_0$ and $\epsilon$ are some positive constants, $c_n$ is a vector such that  $\lim\sup ||c_n|| \to 0$ and $d_n$ is a real number such that $\limsup d_n \to M$ for some constant $M$.
    \item[(C2)] For any $\varepsilon > 0$,
    \[
    \liminf_{n \to \infty} \inf_{\theta \in B_\varepsilon(\theta_n)^c} \{ n^{-1}  (f_n(\theta_n)- f_n(\theta) )\}> 0.
    \]
    \item[(C3)] The prior density $\pi(\theta)$  is continuous at the true value $\theta_0$ and $\pi(\theta_0) > 0$. 
\end{enumerate}


\vspace{-1ex}
The above conditions are introduced to show that the distribution of \(\sqrt{n}(\theta - \theta_n) \mid \bfx_n \) converges to a normal distribution. The density of \( \sqrt{n}(\theta - \theta_n) \mid \bfx_n \) is given by
\[
q_n(x) = \exp\left(f_n\left(\theta_n + \frac{x}{\sqrt{n}}\right) - f_n(\theta_n)\right) \pi\left(\theta_n + \frac{x}{\sqrt{n}}\right) \cdot \frac{\exp(f_n(\theta_n))}{n^{p/2} z_n},
\]
where \( z_n \) is a normalizing constant. The goal is to show that this density converges to that of a normal distribution. 

Condition (C1) ensures that \( f_n(\theta_n + x/\sqrt{n}) - f_n(\theta_n) \) admits a quadratic expansion in \( \theta \), along with control over the remainder term. Equation \eqref{eq:rn1} in Condition (C1) guarantees that the remainder term vanishes pointwise. Equation \eqref{eq:rn2} in Condition (C1), together with Condition (C2), provides an upper bound necessary for applying the dominated convergence theorem: the former gives a bound in a neighborhood of \( \theta_n \), while the latter controls the behavior outside that neighborhood.
Condition (C3) is required for the factor $\pi\left(\theta_n + \frac{x}{\sqrt{n}}\right)$ to converges to $\pi(\theta_0)$.

Theorem \ref{thm:asymptotic} states that the density $q_n$ of $\sqrt{n}(\theta-\theta_n)\mid \bbX$ converges to $\mathcal{N}(0, H_0^{-1})$ in total variation norm, where $\mathcal{N}(0, H_0^{-1})$ is the multivariate normal distribution with mean $0$ and variance $H_0^{-1}$. The structure of the proof closely follows that of Theorem 4 in \citeasnoun{miller2021asymptotic}. However, we provide a detailed argument in the supplementary material because our assumptions are slightly weaker than those imposed in their setting.


\begin{theorem}\label{thm:asymptotic}
Suppose conditions (C1)-(C3) hold and $z_n$ is finite. 
Then, we have 
\[
\int_{\mathbb{R}^p} \left| q_n(x) - \mathcal{N}(x \mid 0, H_0^{-1}) \right| dx \xrightarrow{n \to \infty} 0,
\]
with probability $1$.
\end{theorem}

The above theorem holds for the basic model if we define
\bea 
f_n(\theta):= \log \frac{h(\theta:g^*(F),\alpha+ nF_n)}{h(\theta: g^*(F),\alpha)}. 
\eea

\section{Applications}\label{appl}
In this section, we present two applications of the proposed cDP prior: quantile estimation and moment estimation, while an application for the regression model is provided in the supplementary note. All datasets and R scripts used in this paper are available at \url{https://github.com/statjs/cDPprior}.
\vspace{-3ex}
\subsection{Quantile estimation}
\subsubsection{Posterior of quantiles}
Suppose we observe a random sample from $F$, a distribution on $\mathbb{R}$, and we are interested in estimation of $k \in \bbN$ quantiles $\theta = (\theta_1, \ldots, \theta_k)$  of $F$, where  $\theta_i$ is the $p_i$ quantile of $F$ with $0 < p_1 < \ldots p_k < 1$ and $p = (p_1,p_2, \ldots, p_k)$. This quantile estimation problem can also be formulated as a special case of the functional condition model \eqref{eq:model}. 
Let $\Delta p =  (\Delta p_1, \ldots, \Delta p_{k+1})^T$ with 
$\Delta p_1 = p_1, \Delta p_2 = p_2 - p_1, \ldots, \Delta p_{k+1} = 1 - p_{k}$ and let $S_1(\theta) = (-\infty, \theta_1], S_2(\theta) = (\theta_1, \theta_2], \ldots, S_{k+1} = (\theta_k, \infty)$. 
Define the functional condition as 
\[ g(F, \theta)  = \begin{pmatrix} F(S_1(\theta))  \\ F(S_2(\theta) )  \\ \vdots  \\
		F(S_{k+1}(\theta))  
	\end{pmatrix} = \Delta p. \]
We place the prior on $(F, \theta)$ as
\[ (F, \theta) \sim \mathcal{D}_\alpha(dF | g(F, \theta) = \Delta p) \pi(\theta) d\theta,\]
where $\alpha$ is a finite nonnull measure on $\mathbb{R}$ and $\pi(\theta)$ is a prior density for $\theta$. 
In summary, the quantile model is as follows: 
\begin{equation}\label{eq:quantile_model}
\begin{aligned}
    \theta & \sim \pi(\theta) d\theta \\
    F & \sim \mathcal{D}_\alpha(dF | g(F, \theta) = \Delta p) \\
    y_1, y_2, \ldots, y_n \mid F &\iid F.
\end{aligned}
\end{equation}
Note that when $F \sim \mathcal{D}_\alpha$, 
\[ g(F, \theta) \sim Dir(\alpha(S_1(\theta)), \alpha(S_2(\theta)), \ldots, \alpha(S_{k+1}(\theta)) ). \]
Applying Theorem \ref{th:post_cDP}, we obtain the posterior 
\begin{align}
    \pi(dF, d\theta | \bfy_n) & \propto 
    \mathcal{D}_{\alpha + n F_n} (dF | g(F, \theta) = \Delta p ) 
    \frac{h(\Delta p: g(F, \theta), \alpha + n F_n)} {h(\Delta p: g(F, \theta), \alpha)} \pi(\theta) d\theta \\
    & = \mathcal{D}_{\alpha + n F_n} (dF | g(F, \theta) = \Delta p ) \frac{[\alpha(\mathbb{R})]_n}{\prod_{i=1}^{k+1} [\alpha(S_i(\theta))]_{nF_n(S_i(\theta))} }  \prod_{i=1}^{k+1} \Delta p_i^{ n F_n(S_i(\theta)) }\pi(\theta) d\theta, \nonumber
\end{align}
where $F_n = \frac{1}{n} \sum_{i=1}^n \delta_{y_i}$. 
Thus, the conditional posterior $F$ given $\theta$ and $\bfy_n$  is 
\[ F | \theta, \bfy_n \sim  \mathcal{D}_{\alpha + n F_n} (dF | g(F, \theta) = \Delta p ) \]
and the marginal posterior of $\theta$ is 
\bean\label{eq:quantilemarginal} 
 \pi(d\theta | \bfy_n) \propto \frac{[\alpha(\mathbb{R})]_n}{\prod_{i=1}^{k+1} [\alpha(S_i(\theta))]_{nF_n(S_i(\theta))} }  \prod_{i=1}^{k+1} \Delta p_i^{ n F_n(S_i(\theta)) }\pi(\theta) d\theta. 
 \eean 
Note that the marginal likelihood of $\theta$ is 
\begin{equation} \label{eq:marginal_likelihood}
    L_{cDP} (\theta) = \frac{[\alpha(\mathbb{R})]_n}{\prod_{i=1}^{k+1} [\alpha(S_i(\theta))]_{nF_n(S_i(\theta))} }  \prod_{i=1}^{k+1} \Delta p_i^{ n F_n(S_i(\theta)) }.
\end{equation}
The marginal likelihood of $\theta$ is reminiscent of Jeffreys likelihood for quantiles, 
\begin{equation}
    L_J(\theta) = \frac{\Gamma(n)}{\prod_{i=1}^{k+1} (nF_n(S_i(\theta)))!}   \prod_{i=1}^{k+1} \Delta p_i^{ n F_n(S_i(\theta)) }.
\end{equation}
See \citeasnoun{jeffreys1998theory} and \citeasnoun{lavine1995approximate}.  In fact, the Jeffreys likelihood can be obtained as the limit of the marginal likelihood $L_{cDP}(\theta)$ of $\theta$ by letting $\alpha(\bbR)$ and $\alpha(S_i(\theta))$'s tend to $1$; that is, 
$\lim_{\alpha(\bbR) \to 1, \alpha(S_i(\theta)) \to 1, i=1,2, \ldots, k+1} L_{cDP}(\theta) = L_J(\theta). $
This appears to be the first direct justification for the Jeffreys likelihood. 
A peculiar feature of the marginal likelihood, which is shared by the Jeffreys likelihood,  is that it does not vanish as $\theta \to \pm \infty$. When $k = 1$, 
\begin{align*}
    \lim_{\theta \to \infty} L_{cDP}(\theta)  = p^n,~\lim_{\theta \to -\infty} L_{cDP}(\theta)  = (1-p)^n. 
\end{align*}
Thus, the uniform prior on $\mathbb{R}$ of $\theta$ leads an improper posterior. Therefore, it is necessary to use a proper prior for $\theta$.

We have examined the asymptotic behavior of the posterior distribution under two regimes. First, as the prior mass \( A \to 0 \), the marginal posterior distribution converges to a limiting form proportional to
\[
\pi(\theta \mid \mathbf{y}_n) \lra 
  \frac{1}{z_0 \prod_{j \in \bm{J}} \bar{\alpha}(S_j(\theta)) \Gamma(n F_n(S_j(\theta)))}
  \prod_{j=1}^{k+1} (\Delta p_j)^{n F_n(S_j(\theta))} \pi(\theta), \quad \forall \theta,
\]
where \( z_0 \) is a finite normalizing constant. In particular, if the base measure \( \bar{\alpha} \) is chosen to be the empirical distribution \( F_n \), this limiting posterior becomes proportional to the Jeffreys likelihood multiplied by the prior \( \pi(\theta) \), thereby providing a second justification for the Jeffreys likelihood in quantile inference. Second, as the sample size \( n \to \infty \), the posterior distribution of \( \sqrt{n}(\theta - q_n) \), centered at the sample quantiles $q_n$, converges in total variation to a multivariate normal distribution. The derivation of this result is provided in the supplementary material.

\subsubsection{Posterior computation}

Given a proper prior, the marginal posterior distribution of $\theta$ can be explored using Markov chain Monte Carlo (MCMC) methods. In particular, we adopt the slice sampling algorithm, which is well-suited for target distributions with intractable normalization constants and nonstandard forms. The following algorithm outlines the procedure used to draw posterior samples of $\theta$ from $\pi(d\theta \mid \mathbf{y}_n)$.

\begin{algorithm}[H]
\caption{Slice sampler for the marginal posterior distribution of $\theta$}
\KwIn{Initial value $\theta^{(0)}$, number of iterations $T$}
\KwOut{Samples $\{\theta^{(t)}\}_{t=1}^T$ from the posterior $\pi(\theta \mid \mathbf{y}_n)$}

\For{$t = 1$ to $T$}{
  Draw $u \sim \text{Uniform}(0, \pi(\theta^{(t-1)} \mid \mathbf{y}_n))$\;
  Find interval $(L, R)$ such that $\pi(\theta \mid \mathbf{y}_n) > u$\;
  Sample $\theta^{(t)} \sim \text{Uniform}(\{ \theta : \pi(\theta \mid \mathbf{y}_n) > u \} \cap (L, R))$\;
}
\Return{$\{\theta^{(t)}\}_{t=1}^T$}
\label{alg:quantile_slice}
\end{algorithm}

\subsubsection{Data analysis}\label{sec:quantile_application}
In this study, we present a simple example of estimating quantiles for the maximum sustained wind speeds of hurricanes. Quantile estimation is useful because it helps describe both how likely strong winds are and how intense they can be. While analyses based on averages may miss the effects of rare but powerful hurricanes, quantiles provide more detailed information about the extreme values in the data. This helps us better assess the risk of severe hurricanes and understand how they relate to climate factors.

The dataset used in this example is available at \texttt{https://myweb.fsu.edu/jelsner/temp/}. The sample includes 2,297 observations.
To estimate the quantiles of maximum sustained wind speeds from hurricanes, we denote the observed data as $\mathbf{y}_n = (y_1, \ldots, y_n)$ and assume it follows an unknown distribution $F$. We place a prior on both the distribution $F$ and the quantile $\theta$ as follows: $\pi(dF, d\theta) = \mathcal{D}_\alpha(dF \mid F(\theta) = p)\pi(\theta)d\theta$, 
where $p$ represents the quantile level of interest, such as $0.95$ or $0.99$, corresponding to extreme wind speeds. The base measure is defined as $\alpha(\cdot) := M F_0(\cdot)$, with $M = 1$ and $F_0(\cdot)$ taken as a Cauchy distribution $\mathcal{C}(\cdot \mid \mu_0, \sigma_0)$, and $\pi(\theta) = \mathcal{C}(\theta \mid 0, \tau)$.

Under this formulation, the marginal posterior distribution of the quantile $\theta$ is given by
\begin{equation}\label{eq:post_quantile}
	\pi(d\theta \mid \mathbf{y}_n) \propto \frac{1}{[\alpha(S_\theta)]_{n(S_\theta)}[\alpha(S_\theta^c)]_{n(S_\theta^c)}}p^{n(S_\theta)}(1 - p)^{n(S^c_\theta)}\mathcal{C}(\theta \mid 0, \tau)d\theta,
\end{equation}
where $n(S_\theta) = \sum_{i=1}^n I(y_i \le \theta)$, $n(S_\theta^c) = \sum_{i=1}^n I(y_i > \theta)$, $\alpha(S_\theta) = \int_{-\infty}^\theta \mathcal{C}(y \mid \mu_0, \sigma_0)dy$, and hyper-parameters $\mu_0 = \text{median}(\mathbf{y}_n)$, $\sigma_0 = \text{IQR}(\mathbf{y}_n)/2$, and $\tau = 10 \times \text{sd}(\hat{q}_p) \times 2$ are chosen based on the observed data.

\begin{figure}[!ht]
	\centering
	\includegraphics[width=0.6\textwidth]{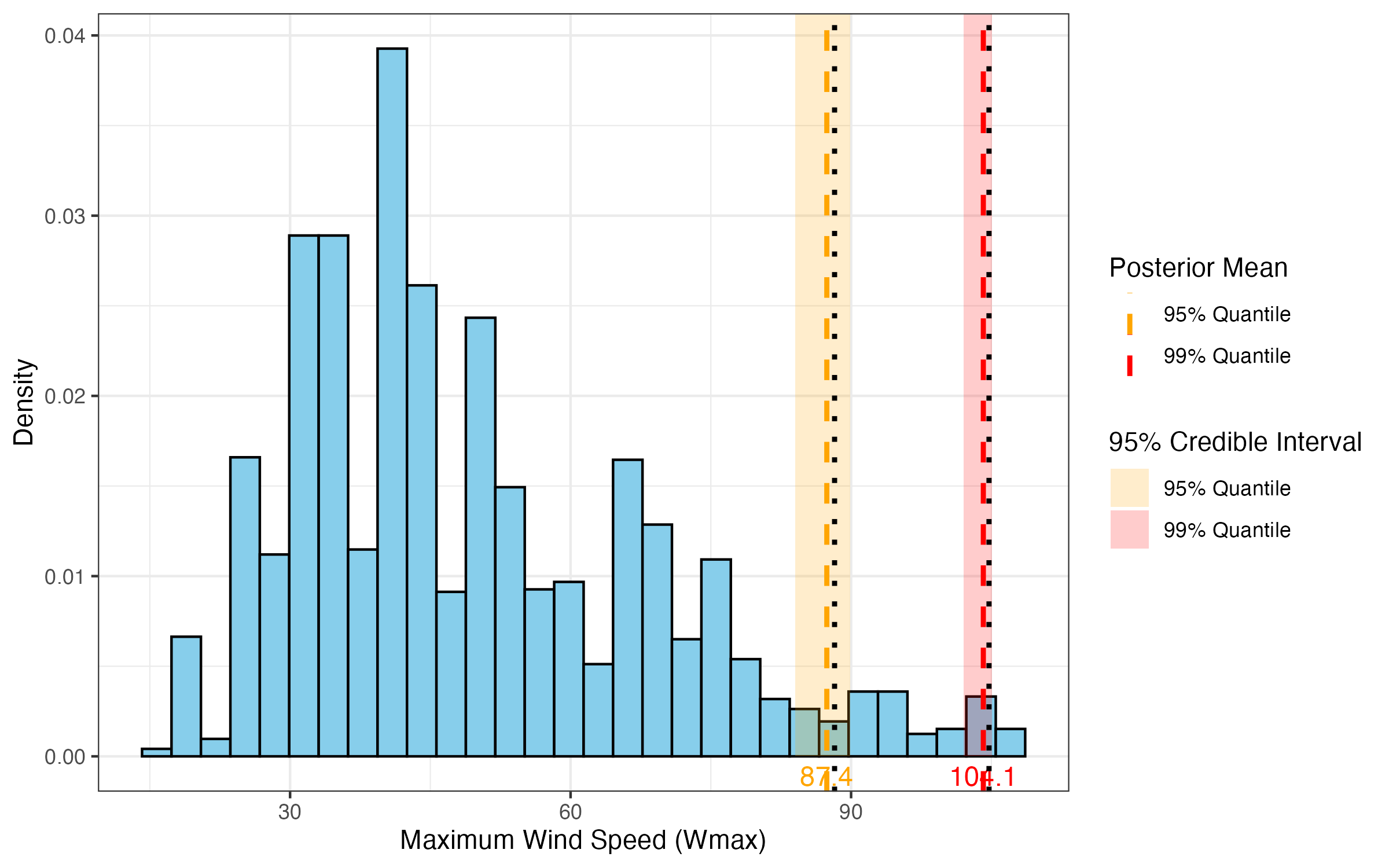} 
	\caption{Posterior summaries of the 95\% and 99\% quantiles of maximum wind speed in 2006. Black dashed lines represent sample quantile estimates.}
	\label{fg:post_wmax}
\end{figure}

Table~\ref{tab:quantile_ci} and Figure~\ref{fg:post_wmax} summarize the estimated quantiles. For the frequentist approach, the 95\% confidence intervals were computed using the bootstrap method. For our approach, we used 10,000 posterior samples, discarding the first 5,000 as burn-in and using the remaining 5,000 for estimation. In the figure, the orange and red vertical lines indicate the estimated 95\% and 99\% quantiles, respectively. For the our method, the sample quantile estimates are shown in black. 
\begin{table}[!ht]
	\centering
	\caption{Quantile estimates and 95\% credible intervals for maximum wind speed in 2006.}
	\begin{tabular}{ll|rrr}
		\hline\hline
		Quantile & Method & Estimate & 95\% CI Lower & 95\% CI Upper \\
		\hline
		95\%     & Sample Quantile &  88.2             & 84.1                   & 90.0 \\ 
		& cDP prior    &  87.4             & 84.0                   & 89.9 \\ 
		\hline
		99\%     & Sample Quantile & 105.0             & 103.2                  & 105.0 \\
		& cDP prior    & 104.0             & 102.0                  & 105.0 \\
		\hline\hline
	\end{tabular}
	\label{tab:quantile_ci}
\end{table}

From the results, for the 95\% quantile, the sample quantile estimate is 88.2, with a 95\% confidence interval of [84.1, 90.0]. In comparison, the cDP-based estimate is slightly lower at 87.4, with a similar credible interval of [84.0, 89.9]. 
For the 99\% quantile, the sample estimate is 105.0, with a confidence interval of [103.2, 105.0], whereas the cDP prior yields an estimate of 104.0 and a slightly wider interval of [102.0, 105.0].
Overall, the cDP-based estimates are close to the frequentist sample quantiles, but tend to be slightly more conservative, especially in the upper quantile.

\subsection{Moment estimation}
The moment model can be formulated as a special case of the standard model \eqref{eq:standard_model}.  
We observed independent copies of $y$, denoted by $y_1, \ldots, y_n \iid F$, where $F$ is a distribution function on $\bbR$.
The goal is to estimate the four moments $\theta = (\mu, \sigma, \gamma, \kappa)$ of $F$ and the functional condition is defined as 
\[ g^*(F)  = \begin{pmatrix} \E(y) \\ sd(y) \\ \E \Big( \frac{y-\mu}{\sigma}\Big)^3 \\ \E \Big( \frac{y-\mu}{\sigma}\Big)^4 
\end{pmatrix} = \theta.   \]
We place the prior on $(F, \theta)$ as $(F, \theta) \sim \mathcal{D}_\alpha(dF| g^*(F) = \theta) \pi(\theta) d\theta$,
where $\alpha$ is a finite nonnull measure on $\mathbb{R}$ and $\pi(\theta)$ is a prior density for $\theta$. The moment model is specified as follows: 
\begin{equation}\label{eq:moment_model}
	\begin{aligned}
		\theta & \sim \pi(\theta) d\theta \\
		F & \sim \mathcal{D}_\alpha (dF | g^*(F) = \theta) \\
		y_1, y_2, \ldots, y_n & \iid F. 
	\end{aligned}
\end{equation}
The posterior distribution is given as follows. 
\begin{equation} \label{eq:post_moment_model}
	\begin{aligned}
		\pi(dF | \theta, \bfy_n) & = \mathcal{D}_{\alpha + nF_n} (dF| g^*(F) = \theta) \\
		\pi(d\theta | \bfy_n) & \propto \frac{h(\theta: g^*(F), \alpha + nF_n)}{h(\theta: g^*(F), \alpha)} \pi(\theta) d\theta.
	\end{aligned}
\end{equation}

Note that if we use $h(\theta: g^*(F), \alpha)$ as the prior for $\theta$, then the posterior of $\theta$ is given by $\pi(d\theta | \bfy_n) \propto h(\theta: g^*(F), \alpha + nF_n) d\theta$,
and by letting $A \to 0$, we obtain the Bayesian bootstrap posterior    $\pi^{BB}(d\theta|\bfy_n) \propto  h(\theta: g^*(F), nF_n) d\theta$.

\subsubsection{Posterior computation}\label{sec:post_comp_moment}
Given a base measure $\alpha(\cdot)$ and a proper prior $\pi(\theta)$, we can generate posterior samples of $\theta$ by drawing from constrained distributions derived from the Bayesian bootstrap \cite{rubin1981bayesian}, and applying importance weighting to adjust for the prior. The following algorithm outlines the procedure.

\begin{algorithm}[H]
	\caption{Importance-reweighted sampling for the marginal posterior distribution of $\theta$}
	\KwIn{Number of samples $T$}
	\KwOut{Posterior samples $\{\theta^{(t)}\}_{t=1}^T$ from $\pi(\theta \mid \mathbf{y}_n)$}
	\For{$t = 1$ \KwTo $T$}{
		Sample $F_{\text{prior}}^{(t)} \sim \mathcal{D}_{\alpha}$ via stick-breaking with $N$ atoms from $F_0$\;
		Sample $F_{\text{data}}^{(t)} \sim \mathcal{D}_{nF_n}$ using the Bayesian bootstrap on $\mathbf{y}_n$\;
		Draw mixing weight $\lambda^{(t)} \sim \text{Beta}(n, M)$\;
		Form $F^{(t)} = \lambda^{(t)} F_{\text{data}}^{(t)} + (1 - \lambda^{(t)}) F_{\text{prior}}^{(t)}$\;
		Solve for $\theta^{(t)}$ such that $g^*(F^{(t)}) = \theta^{(t)}$\;
		Compute weight $w^{(t)} \propto \pi(\theta^{(t)}) / h(\theta^{(t)}; \alpha)$\;
	}
	Normalize weights $\{w^{(t)}\}$ and resample $\theta^{(t)}$ with replacement\;
	\Return{$\{\theta^{(t)}\}_{t=1}^T$}
	\label{alg:importance_sampling_moments}
\end{algorithm}

\subsubsection{Data analysis}
\begin{figure}[!ht]
	\centering
	\includegraphics[width=0.7\textwidth]{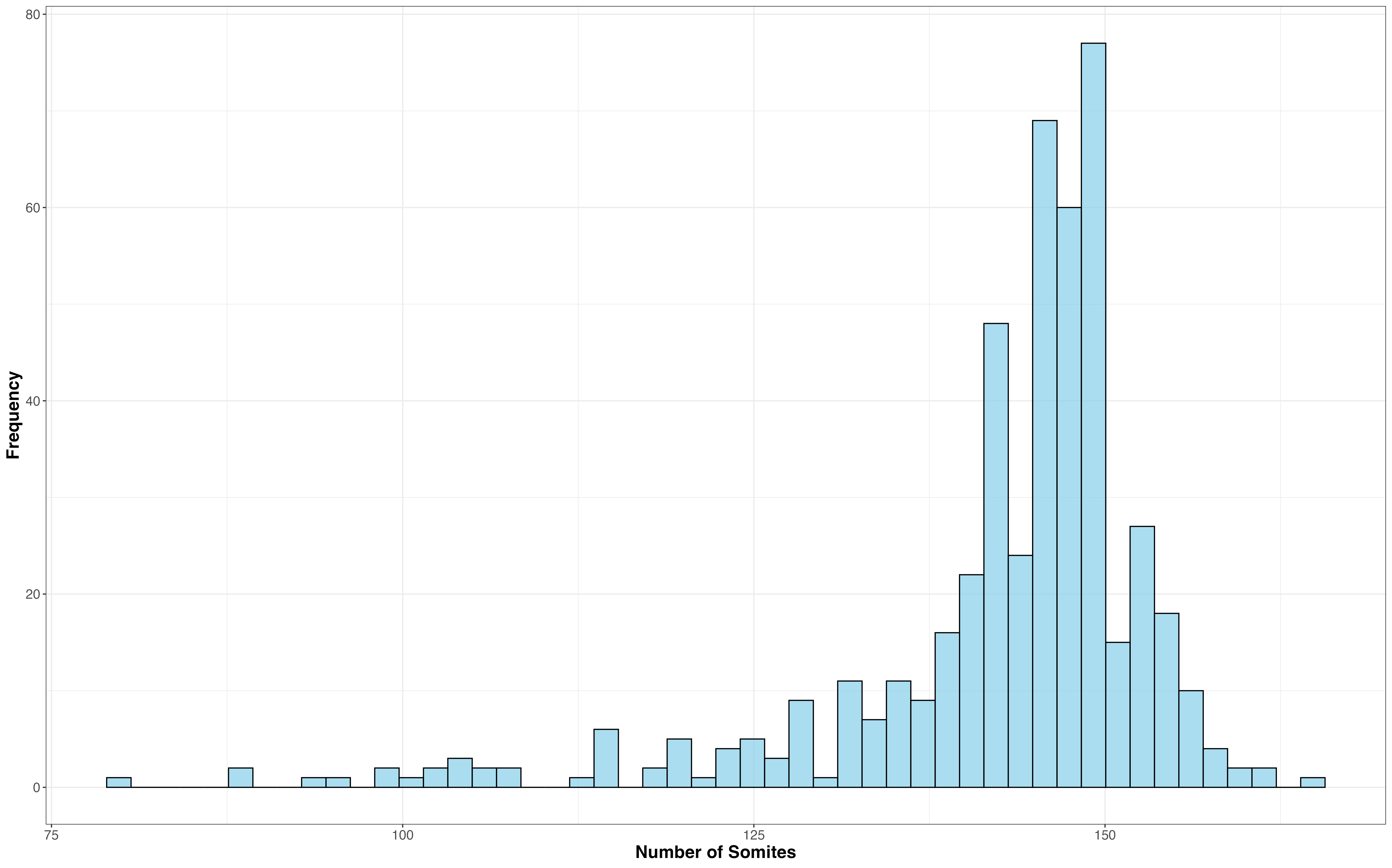}
	\caption{Histogram of the number of somites observed in 487 earthworms.}
	\label{fg:hist_somite}
\end{figure}

As an example, we consider a classical dataset originally analyzed by \citeasnoun{pearl1905earthworm}, which reports the number of somites, the segmented units of the body, in individual earthworms. The somite count serves as a biologically meaningful measure of morphological complexity and developmental variation among specimens. The dataset consists of 487 individual observations of earthworms. Each observation records the number of somites along with its corresponding frequency. The histogram (Figure \ref{fg:hist_somite}) displays the empirical distribution of the somite counts. The distribution is approximately unimodal and slightly right-skewed.

To estimate the four moments $\theta = (\mu, \sigma, \gamma, \kappa)$, we assume a base measure as \(\alpha(\cdot) = t_\nu(\cdot \mid \mu_0, \sigma_0)\), where \(\mu_0 = \text{median}(\mathbf{y}_n)\) and \(\sigma_0 = \text{IQR}(\mathbf{y}_n) / 1.454\), with the degrees of freedom fixed at \(\nu = 5\) to ensure the existence of the fourth moment. The prior distributions are specified as follows: \(\mu \sim \mathcal{C}(m_1, \tau_1)\), \(\sigma \sim \mathcal{G}a(a_2, b_2)\), \(\gamma \sim \mathcal{C}(m_3, \tau_3)\), and \(\kappa \sim \mathcal{G}a(a_4, b_4)\). Under this prior specification, the posterior is proportional to
\begin{equation*}
	\pi(d\theta | \bfy_n )  \propto      \frac{h(\theta: g^*(F), \alpha + n F_{n} )} { h(\theta: g^*(F), \alpha)} \mathcal{C}(\mu : m_1, \tau_1)\mathcal{G}a(\sigma : a_2, b_2)\mathcal{C}(\gamma : m_3, \tau_3)\mathcal{G}a(\kappa : a_4, b_4)d\theta, 
\end{equation*}
where the hyperparameters, $m_1, \tau_1, a_2, b_2, m_3, \tau_3, a_4$, and $a_5$, are empirically specified based on the observed data, similarly to the approach in Section \ref{sec:quantile_application}.

\begin{figure}[!ht]
	\centering
	\includegraphics[width=0.7\textwidth]{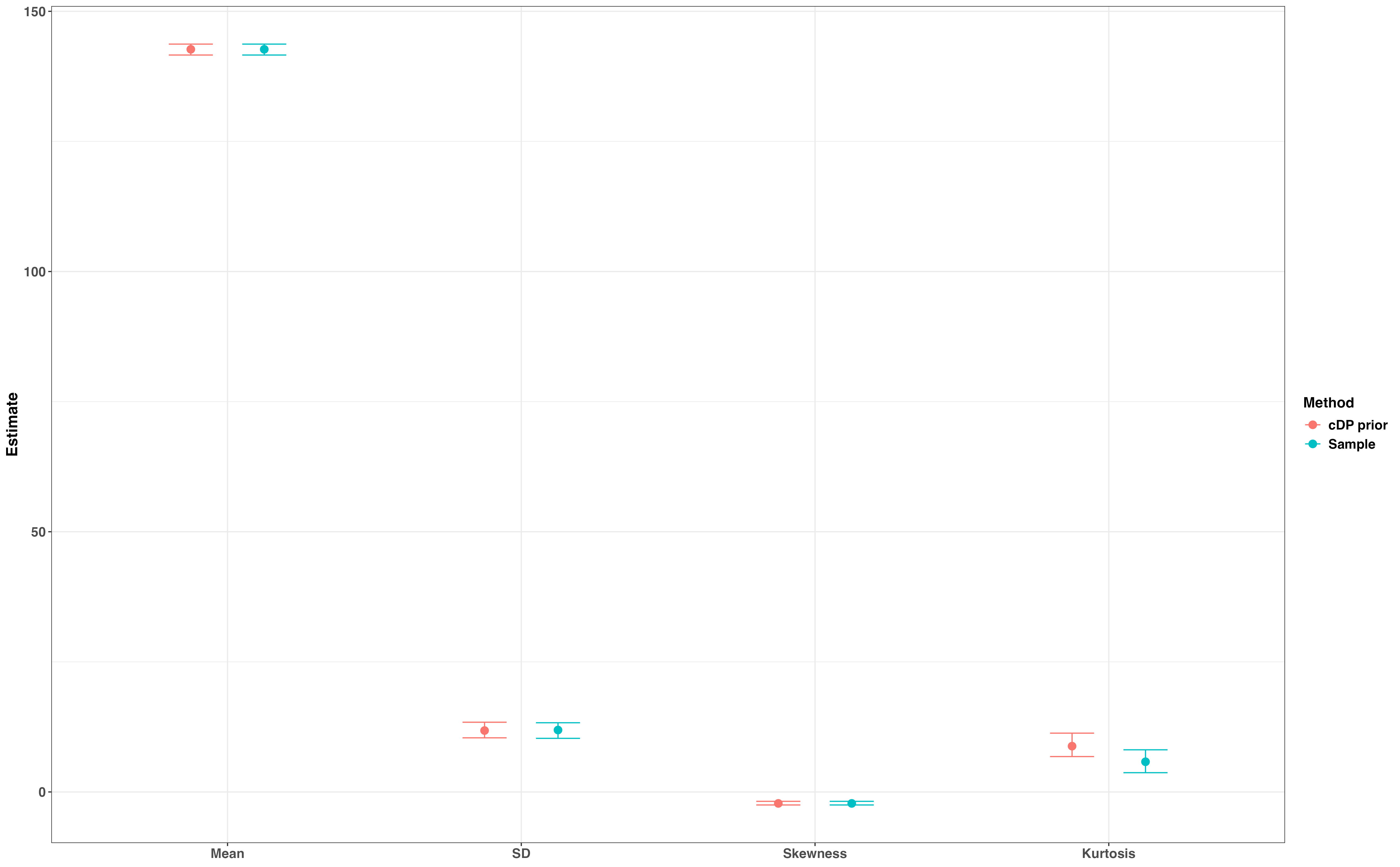} 
	\caption{Posterior summaries of the four moments estimates of the number of somites.}
	\label{fg:post_moment}
\end{figure}

Table~\ref{tab:moment_ci} summarizes the estimated values of the first four moments of the somite count distribution in earthworms, along with their associated interval estimates, and these results are visualized in Figure~\ref{fg:post_moment}. For comparison, we report both sample-based estimates, obtained via nonparametric bootstrap with 1,000 resamples, and posterior estimates based on 1,000 posterior draws under the proposed model, with trace plot diagnostics indicating good convergence for all parameters of interest.
\begin{table}[!ht]
	\centering
	\caption{Moment estimates and 95\% credible intervals for earthworms.}
	\begin{tabular}{ll|rrr}
		\hline\hline
		Moment & Method & Estimate & 95\% CI Lower & 95\% CI Upper \\
		\hline
		$\mu$     & Sample mean &  142.7             & 141.6                   & 143.7 \\ 
		& cDP prior    &  142.7             & 141.6                   & 143.7 \\ 
		\hline
		$\sigma$     & Sample std & 11.9             & 10.3                  & 13.3 \\
		& cDP prior    & 11.8            & 10.4                  & 13.4 \\
		\hline
		$\gamma$     & Sample skewness & -- 2.2             & -- 2.5                  & -- 1.8 \\
		& cDP prior    & -- 2.2            & -- 2.5                 & -- 1.8 \\
		\hline
		$\kappa$     & Sample kurtosis & 5.8             & 3.7                  & 8.1 \\
		& cDP prior    & 8.8             & 6.8                  & 11.3 \\
		\hline\hline
	\end{tabular}
	\label{tab:moment_ci}
\end{table}

Both methods yield nearly identical estimates of the mean somite count (\(\mu = 142.7\)), with matching uncertainty intervals \([141.6,\ 143.7]\). Similarly, the posterior standard deviation is very close across methods (\(\sigma = 11.9\) for the sample vs.\ 11.8 for the cDP prior), with overlapping intervals. The skewness estimates are also identical (\(\gamma = -2.2\), 95\% CI \([-2.5,\ -1.8]\)), suggesting consistent asymmetry detection by both approaches. In contrast, the cDP prior yields a substantially higher kurtosis estimate (\(\kappa = 8.8\) vs.\ 5.8), indicating heavier tails. These results suggest that the cDP prior may better capture the presence of potential outliers, such as specimens with unusually low somite counts, that are less clearly reflected in sample-based estimation alone.

\section{Discussion}\label{disc}

This study has developed a non-parametric Bayesian methodology for models in which the parameter of interest is defined as the solution to a functional equation involving an unknown distribution. We introduced the functional condition model and formulated the conditional Dirichlet process (cDP) as a prior tailored to this structure. This framework allows priors to be placed directly on parameters implicitly defined by a functional of the distribution and facilitates posterior inference by deriving the marginal posterior distribution of the parameter and the conditional posterior distribution of the distribution. We also proposed a practical approach to posterior computation based on this formulation. Theoretical results were established for the limiting behaviour of the marginal posterior of $\theta$, including its asymptotic form as the prior mass tends to zero, and its asymptotic normality as the sample size increases. The proposed methodology was applied to quantile and moment models. A natural direction for future work is to extend this framework to regression models, such as quantile and instrumental regression, where the regression coefficients are defined via a functional condition.

\bibliographystyle{dcu}           
\bibliography{cdp}    

@article{pearl1905earthworm,
	author = {Pearl, Raymond},
	title = {On the Distribution of the Number of Somites in Earthworms},
	journal = {Biometrika},
	year = {1905},
	volume = {4},
	number = {1},
	pages = {40--52}
}

@article{miller2021asymptotic,
  title={Asymptotic normality, concentration, and coverage of generalized posteriors},
  author={Miller, Jeffrey W},
  journal={Journal of Machine Learning Research},
  volume={22},
  number={168},
  pages={1--53},
  year={2021}
}

@article{robbins1955remark,
  title={A remark on Stirling's formula},
  author={Robbins, Herbert},
  journal={The American mathematical monthly},
  volume={62},
  number={1},
  pages={26--29},
  year={1955},
  publisher={JSTOR}
}

@article{alam2025dpglm,
  title={DPGLM: A Semiparametric Bayesian GLM with Inhomogeneous Normalized Random Measures},
  author={Alam, Entejar and Rathouz, Paul J and Mueller, Peter},
  journal={arXiv:2502.17827},
  year={2025}
}

@article{chang1997conditioning,
  title={Conditioning as disintegration},
  author={Chang, Joseph T and Pollard, David},
  journal={Statistica Neerlandica},
  volume={51},
  number={3},
  pages={287--317},
  year={1997},
  publisher={Wiley Online Library}
}

@article{chaudhuri2008generalized,
  title={Generalized linear models incorporating population level information: an empirical-likelihood-based approach},
  author={Chaudhuri, Sanjay and Handcock, Mark S and Rendall, Michael S},
  journal={Journal of the Royal Statistical Society Series B: Statistical Methodology},
  volume={70},
  number={2},
  pages={311--328},
  year={2008},
  publisher={Oxford University Press}
}

@article{chaudhuri2017hamiltonian,
  title={Hamiltonian Monte Carlo sampling in Bayesian empirical likelihood computation},
  author={Chaudhuri, Sanjay and Mondal, Debashis and Yin, Teng},
  journal={Journal of the Royal Statistical Society Series B: Statistical Methodology},
  volume={79},
  number={1},
  pages={293--320},
  year={2017},
  publisher={Oxford University Press}
}

@article{chib2018bayesian,
  title={Bayesian estimation and comparison of moment condition models},
  author={Chib, Siddhartha and Shin, Minchul and Simoni, Anna},
  journal={Journal of the American Statistical Association},
  volume={113},
  number={524},
  pages={1656--1668},
  year={2018},
  publisher={Taylor \& Francis}
}

@article{chib2022bayesian,
  title={Bayesian estimation and comparison of conditional moment models},
  author={Chib, Siddhartha and Shin, Minchul and Simoni, Anna},
  journal={Journal of the Royal Statistical Society Series B: Statistical Methodology},
  volume={84},
  number={3},
  pages={740--764},
  year={2022},
  publisher={Oxford University Press}
}

@article{ferguson1973bayesian,
  title={A Bayesian analysis of some nonparametric problems},
  author={Ferguson, Thomas S},
  journal={The annals of statistics},
  pages={209--230},
  year={1973},
  publisher={JSTOR}
}

@book{ghosal2017fundamentals,
  title={Fundamentals of nonparametric Bayesian inference},
  author={Ghosal,  Subhashis and Van der Vaart,  Aad},
  year={2017},
  publisher={Cambridge University Press}
}

@book{jeffreys1998theory,
  title={The theory of probability},
  author={Jeffreys, Harold},
  year={1998},
  publisher={OuP Oxford}
}

@article{lavine1995approximate,
  title={On an approximate likelihood for quantiles},
  author={Lavine, M},
  journal={Biometrika},
  volume={82},
  number={1},
  pages={220--222},
  year={1995},
  publisher={Oxford University Press}
}

@article{lazar2003bayesian,
  title={Bayesian empirical likelihood},
  author={Lazar, Nicole A},
  journal={Biometrika},
  volume={90},
  number={2},
  pages={319--326},
  year={2003},
  publisher={Oxford University Press}
}

@book{owen2001empirical,
  title={Empirical likelihood},
  author={Owen, Art B},
  year={2001},
  publisher={Chapman and Hall/CRC}
}

@article{porter2015bayesian,
  title={Bayesian semiparametric hierarchical empirical likelihood spatial models},
  author={Porter, Aaron T and Holan, Scott H and Wikle, Christopher K},
  journal={Journal of Statistical Planning and Inference},
  volume={165},
  pages={78--90},
  year={2015},
  publisher={Elsevier}
}

@article{schennach2005bayesian,
  title={Bayesian exponentially tilted empirical likelihood},
  author={Schennach, Susanne M},
  journal={Biometrika},
  volume={92},
  number={1},
  pages={31--46},
  year={2005},
  publisher={Oxford University Press}
}

@book{tjur1975constructive,
  title={A constructive definition of conditional distributions},
  author={Tjur, Tue},
  year={1975},
  publisher={Institute of Mathematical Statistics, University of Copenhagen}
}

@article{yang2012bayesian,
  title={Bayesian Empirical Likelihood for Quantile Regression},
  author={Yang, Y. and He, X.},
  journal={The Annals of Statistics},
  volume={40},
  pages={1102--1131},
  year={2012}
}

@article{rubin1981bayesian,
  title={The Bayesian bootstrap},
  author={Rubin, Donald B.},
  journal={The Annals of Statistics},
  volume={9},
  number={1},
  pages={130--134},
  year={1981},
  publisher={Institute of Mathematical Statistics}
}

@misc{nhanes2021,
  author       = {{Centers for Disease Control and Prevention (CDC), National Center for Health Statistics (NCHS)}},
  title        = {{National Health and Nutrition Examination Survey: 2021–2023 Data Documentation, Codebook, and Frequencies}},
  year         = {2023},
  howpublished = {\url{https://wwwn.cdc.gov/nchs/nhanes/continuousnhanes/default.aspx?Cycle=2021-2023}},
  note         = {Accessed May 22, 2025}
}

\end{document}


\title{Supplementary Material for "Conditional Dirichlet Processes and Functional Condition Models"}

 \author[1]{Jaeyong Lee}
 \author[2]{Kwangmin Lee}
 \author[3]{Jaegui Lee}
 \author[4]{Seongil Jo}

 \affil[1]{Department of Statistics, Seoul National University}
 \affil[2]{Department of Big Data Convergence, Chonnam National University}
 \affil[3]{SK Innovation}
 \affil[4]{Department of Statistics, Inha University}

\maketitle

\section{Theoretical properties of the posterior in quantile problem}

First, we provide the limiting distribution of the posterior \eqref{eq:quantilemarginal} as $A\to 0$. We need the following assumption. 
\begin{enumerate}
    \item[(QA1)]
The function
        \[
          G(\theta) = 
          \frac{\Gamma(n+1)}{\prod_{j\in \bm{J} }\bar\alpha(S_j(\theta))\Gamma(nF_n(S_j(\theta)))}
          \prod_{j=1}^{k+1}(\Delta p_j)^{nF_n(S_j(\theta))}
          \pi(\theta)
        \]
        is integrable, where $\bm{J} = \{ j\in [k+1] : F_n(S_j(\theta))>0 \} $.
\end{enumerate}

Condition (QA1) is required to ensure condition (B3) for the quantile model \eqref{eq:quantile_model}.
By applying Theorem \ref{thm:asymptoticA} to the quantile model \eqref{eq:quantile_model}, we obtain Theorem~\ref{thm:quantileAasymp} which provides the limiting distribution of the posterior as $A \to 0$.
The proof of Theorem \ref{thm:quantileAasymp} is given in Section \ref{sec:proofquantileasymp}.
\begin{theorem}\label{thm:quantileAasymp}
Suppose the quantile model \eqref{eq:quantile_model} and condition (QA1) hold. Then, as $A\to 0$,
\bean\label{eq:limitingAquantiles} 
  \pi(\theta| \bfy_n) \lra 
  \frac{1}{z_0\prod_{j\in \bm{J}}\bar\alpha(S_j(\theta))\Gamma\bigl(nF_n(S_j(\theta))\bigr)}
  \prod_{j=1}^{k+1}(\Delta p_j)^{nF_n(S_j(\theta))} \pi(\theta), ~\forall\theta
\eean 
where $z_0$ is the finite normalizing constant.
\end{theorem}

Note that even as $A \to 0$, the limiting posterior \eqref{eq:limitingAquantiles} is not free of $\alpha$. However, if we choose   $\alpha := A \times F_n$ with $A > 0$, then Corollary~\ref{corr:quantileAasymp} shows that the limiting posterior is  proportional to
\[
\frac{1}{\prod_{j=1}^{k+1}(nF_n(S_j(\theta)))!}
  \prod_{j=1}^{k+1}(\Delta p_j)^{nF_n(S_j(\theta))}
  \pi(\theta),
\]
which is the Jeffreys likelihood multiplied by the prior $\pi(\theta)$. This provides a second justification for the Jeffreys likelihood.
The proof of Corollary~\ref{corr:quantileAasymp} is given in Section \ref{sec:proofquantileasymp}.

\begin{corollary}\label{corr:quantileAasymp}
Consider the quantile model \eqref{eq:quantile_model} with $\alpha:= A \times F_n$. Then, as $A \to 0$, 
\[
  \pi(\theta| \bfy_n)  \lra
  \frac{1}{z_*\prod_{j=1}^{k+1}(nF_n(S_j(\theta))\bigr)!}
  \prod_{j=1}^{k+1}(\Delta p_j)^{nF_n(S_j(\theta))}
  \pi(\theta),~ \forall\theta
\]
where $z_*$ is a finite normalizing constant.
\end{corollary}

Next, we investigate the limiting distribution of the posterior \eqref{eq:quantilemarginal} as $n \to \infty$. Suppose $y_1, \ldots, y_n$ are generated independently from the true distribution $F_0$ on $\mathbb{R}$,  and let \( \theta_0 = (\theta_{0,1}, \ldots, \theta_{0,k})^\top \) denote the vector of population quantiles. That is, \( \theta_{0,j} \) is the \( p_j \)-quantile of \( F_0 \).
Define \( q_n = (q_{n,1}, \ldots, q_{n,k})^\top \in \mathbb{R}^k \) as the vector of  sample quantiles, where
\[
q_{n,j} := \inf \{ x \in \mathbb{R} : F_n(-\infty,x] \ge p_j \}, \quad j = 1, \ldots, k.
\]
We derive the asymptotic posterior distribution of the quantile \( \theta \), centered at the  sample quantile \( q_n \). To this end, we assume the following conditions. 

\begin{enumerate}
\item[(Qn1)] (Local Lipschitz continuity) There exists a constant \( L > 0 \) such that, for each \( j=1,\dots,k \), $\alpha(x, x'] \le L |x - x'|$ for all $x < x'$  sufficiently close to $\theta_{0,j}$.
\item[(Qn2)] (Non-degeneracy of increments) The increments of \( \alpha \) at \( \theta_0 \) are strictly positive; that is, $\alpha(S_j(\theta_0)) > 0 $ for all $j = 1, \ldots, k+1$.
\item[(Qn3)] (Lower bound relative to \( F_0 \)) There exists a constant \( M' > 0 \) such that, for all \( a < b \),
$\alpha(a,b] \ge M'F_0(a,b]$.
\end{enumerate}

Condition (Qn1) guarantees the quadratic approximation required in condition (C1). Condition (Qn2) ensures non-negligible prior mass around the true parameter \( \theta_0 \). Condition (Qn3) ensures that the prior measure \( \alpha \) adequately supports regions with significant mass under the true distribution \( F_0 \). Specifically, (Qn3) prevents \( \alpha \) from vanishing on regions with substantial \( F_0 \)-mass and ensures its tails are not lighter than those of \( F_0 \). For instance, if the population distribution follows a Cauchy distribution, a base measure with lighter tails than the Cauchy distribution, such as the Gaussian distribution, would fail to satisfy condition (Qn3).

Applying Theorem \ref{thm:asymptotic}, we obtain the following theorem, which provides the asymptotic distribution of \eqref{eq:quantilemarginal} and offers a rigorous asymptotic justification for posterior estimation centered around empirical quantiles.
The proof of Theorem \ref{thm:asympquantiles} is given in Section \ref{sec:proofquantileasymp2}.

\vspace{-2ex}
\begin{theorem}\label{thm:asympquantiles}
Suppose conditions (Qn1), (Qn2), and (Qn3) hold. Assume further that \( F_0 \) is differentiable with derivative \( f_0 \), and that \( f_0 \) is Lipschitz continuous. Then, the posterior distribution of the scaled quantile parameter satisfies $\left\| p_{\sqrt{n}(\theta - q_n)\mid \bfy_n} - \mathcal{N}(0, \Sigma) \right\|_{\mathrm{TV}} \xrightarrow{p} 0$,
where \( p_{\sqrt{n}(\theta - q_n)\mid \bfy_n} \) denotes the posterior density of \( \sqrt{n}(\theta - q_n) \) given the data \( \bfy_n \), and  $\Sigma = \left( T^\top H T \right)^{-1}$.
Here, \( T \in \mathbb{R}^{k \times k} \) is a lower bidiagonal matrix with entries
$
T_{j,j} = f(\theta_{0,j})$, $T_{j,j-1} = -f(\theta_{0,j-1})$,
and $H = \operatorname{diag}\left( \frac{1}{[\Delta p]_1}, \ldots, \frac{1}{[\Delta p]_k} \right) + \frac{1}{1 - p_k} \cdot \bm{1}\bm{1}^\top$.
\end{theorem}

\section{Linear regression model with random covariates}
The goal of the regression model with random covariates is to estimate the parameter vector $\bbeta = (\beta_1, \ldots, \beta_p)$ satisfying the relation
\begin{equation}\label{eq:regmd2}
    y = \bfx'\bbeta + \epsilon, 
\end{equation}
where $(y, \bfx) \sim F$, $F$ is a distribution on $\mathbb{R}^{p+1}$, $y, \epsilon \in \mathbb{R}$ and $\bfx = (x_1, \ldots, x_p)^\top \in \mathbb{R}^{p}$. We observe 
\[
    (y_1, \bfx_1), \ldots, (y_n, \bfx_n) \stackrel{i.i.d.}{\sim} F.
\]
We further assume $\mathbb{E}(\bfx\epsilon)  = 0$. While it is common to additionally assume $\mathbb{E}(\epsilon) = 0$, but we do not require this condition here. Define 
\[ g^*(F) = \E(\bfx \bfx^T) ^{-1} \E (\bfx y). \]

The regression model can be specified as 
\begin{equation} \label{eq:regression_model}
    \begin{aligned}
        \beta & \sim \pi(\beta) d\beta \\
        F & \sim \mathcal{D}_\alpha(dF | g^*(F) =  \beta) \\
        (y_1, \bfx_1), \ldots, (y_n, \bfx_n) & \iid F. 
    \end{aligned}
\end{equation}
The posterior of $(F, \beta)$ is given by 
\begin{equation}
    \begin{aligned}
        \pi(dF| \beta, D) & = \mathcal{D}_{\alpha + nF_n} (dF | g^*(F)=\beta)  \\
        \pi(d\beta | D) & \propto \frac{h(\beta: g^*(F), \alpha + n F_n)}{h(\beta, g^*(F), \alpha)} \pi(\beta) d\beta,
    \end{aligned}
\end{equation}
where $D = \{(y_1, \bfx_1), \ldots, (y_n, \bfx_n)  \}$.

\newpage
\subsection{Posterior computation}
Posterior samples for the regression coefficient $\bbeta$ can be obtained using a procedure similar to the one used for moment estimation given in section \ref{sec:post_comp_moment}. The algorithm below describes the steps for drawing samples from the marginal posterior distribution of $\bbeta$.

\begin{algorithm}[H]
\caption{Posterior sampling algorithm for the regression coefficient $\bbeta$}
\KwIn{Number of posterior samples $T$}
\KwOut{Posterior samples ${\bbeta^{(t)}}_{t=1}^T$ from $\pi(\bbeta \mid D)$}

\For{$t = 1$ \KwTo $T$}{
Sample prior distribution $F_{\text{prior}}^{(t)} \sim \mathcal{D}_\alpha$ via stick-breaking with $N$ atoms drawn from the base measure $F_0$\;
Sample empirical distribution $F_{\text{data}}^{(t)} \sim \mathcal{D}_{nF_n}$ using the Bayesian bootstrap on the observed data ${(y_i, \bfx_i)}_{i=1}^n$\;
Draw a mixing weight $\lambda^{(t)} \sim \text{Beta}(n, \alpha)$\;
Construct the mixture distribution $F^{(t)} = \lambda^{(t)} F_{\text{data}}^{(t)} + (1 - \lambda^{(t)}) F_{\text{prior}}^{(t)}$\;
Compute $\bbeta^{(t)}$ by performing weighted least squares regression under $F^{(t)}$ as
\[
\bbeta^{(t)} = \left( \sum_{i=1}^N w_i^{(t)} \bfx_i^{(t)} \bfx_i^{(t)\top} \right)^{-1} \left( \sum_{i=1}^N w_i^{(t)} \bfx_i^{(t)} y_i^{(t)} \right);
\]
Compute importance weight $w^{(t)} \propto \pi(\bbeta^{(t)}) / h(\bbeta^{(t)}; g^*(F), \alpha)$\;
}

Normalize weights ${w^{(t)}}$ and resample $\bbeta^{(t)}$ with replacement using the normalized weights\;
\Return{${\bbeta^{(t)}}_{t=1}^T$}
\label{alg:posterior_sampling_beta}
\end{algorithm}

\subsection{Data analysis}
As an example, we constructed a regression model with random covariates using data from the National Health and Nutrition Examination Survey (NHANES) 2021–2023, a nationally representative survey administered by the U.S. Centers for Disease Control and Prevention (CDC) \citeasnoun{nhanes2021}. NHANES collects a wide range of health, nutritional, and physiological data through interviews, physical examinations, and laboratory tests. The datasets were obtained from the publicly accessible NHANES website (\texttt{https://wwwn.cdc.gov/nchs/nhanes/}), and included demographic variables, physical measurements, and blood test results. Specifically, the following files were used: DEMO\_L (demographics), BMX\_L (body measurements), BPXO\_L (blood pressure), GLU\_L (fasting glucose), INS\_L (fasting insulin), and FASTQX\_L (fasting time).

We merged these datasets using the common participant identifier \texttt{SEQN}, and variables relevant to glucose metabolism were extracted. The response variable was fasting plasma glucose (\texttt{LBXGLU}), and explanatory variables included body weight (\texttt{BMXWT}), waist circumference (\texttt{BMXWAIST}), systolic blood pressure (\texttt{BPXOSY1}), diastolic blood pressure (\texttt{BPXODI1}), fasting insulin level (\texttt{LBXIN}), age (\texttt{RIDAGEYR}), and hours since last meal (\texttt{PHAFSTHR}). Observations with missing values in any of the selected variables were removed, resulting in a final sample of 3,266 individuals.

\begin{figure}[!ht]
    \centering
    \includegraphics[width=\textwidth]{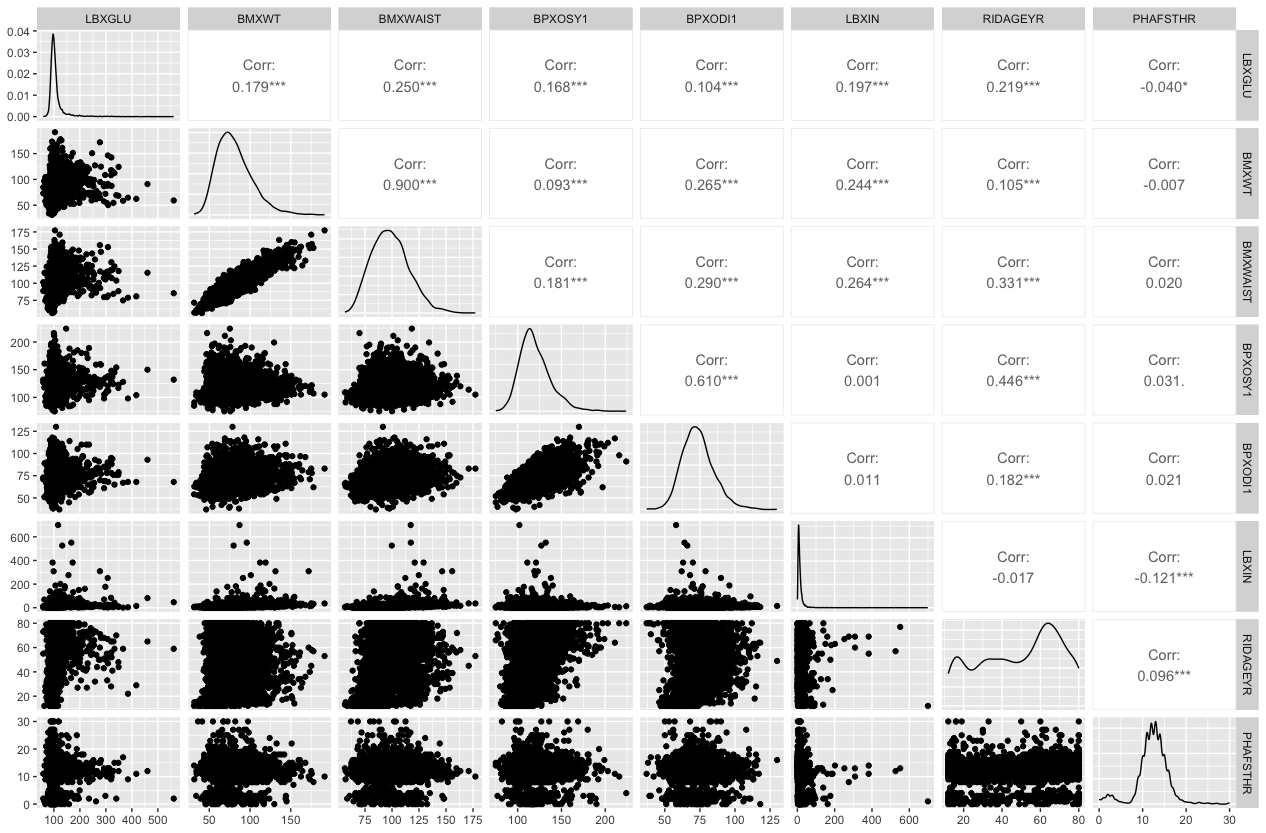}
    \caption{Scatter plot matrix of fasting glucose and explanatory variables.}
    \label{fig:pairs}
\end{figure}

To explore the relationships among the variables, we visualized their pairwise associations using a scatter plot matrix (Figure~\ref{fig:pairs}). Several explanatory variables, such as waist circumference (\texttt{BMXWAIST}) and age (\texttt{RIDAGEYR}), show moderate associations with fasting plasma glucose, while strong correlations are observed between systolic and diastolic blood pressure as well as between body weight and waist circumference.

To estimate the regression coefficient vector $\theta = \bbeta = (\beta_1, \ldots, \beta_p)^\top$, we assume a base measure $\alpha(\cdot) = M \mathcal{C}(\cdot \mid \mathbf{0}, 10^2 \mathbf{I}_{p+1})$, where $M = 1$ is the total mass parameter. We further place independent Cauchy priors on each component of $\bbeta$, that is, $\beta_j \sim \mathcal{N}(\mu_j, \sigma_j)$ for $j = 1, \ldots, p$, where \(\{\mu_j, \sigma_j\}\) are user-specified hyperparameters, selected to be noninformative. Under this specification, the marginal posterior distribution of $\bbeta$ is proportional to
\[
\pi(d\bbeta \mid D) \propto 
\frac{h\left(\beta : g^*(F), \alpha + n F_{n}  \right)}{h\left(\beta : g^*(F), \alpha\right)} 
\prod_{j=1}^p \mathcal{C}(\beta_j \mid \mu_j, \sigma_j)\, d\bbeta.
\]

\begin{figure}[!ht]
    \centering
    \includegraphics[width=\textwidth]{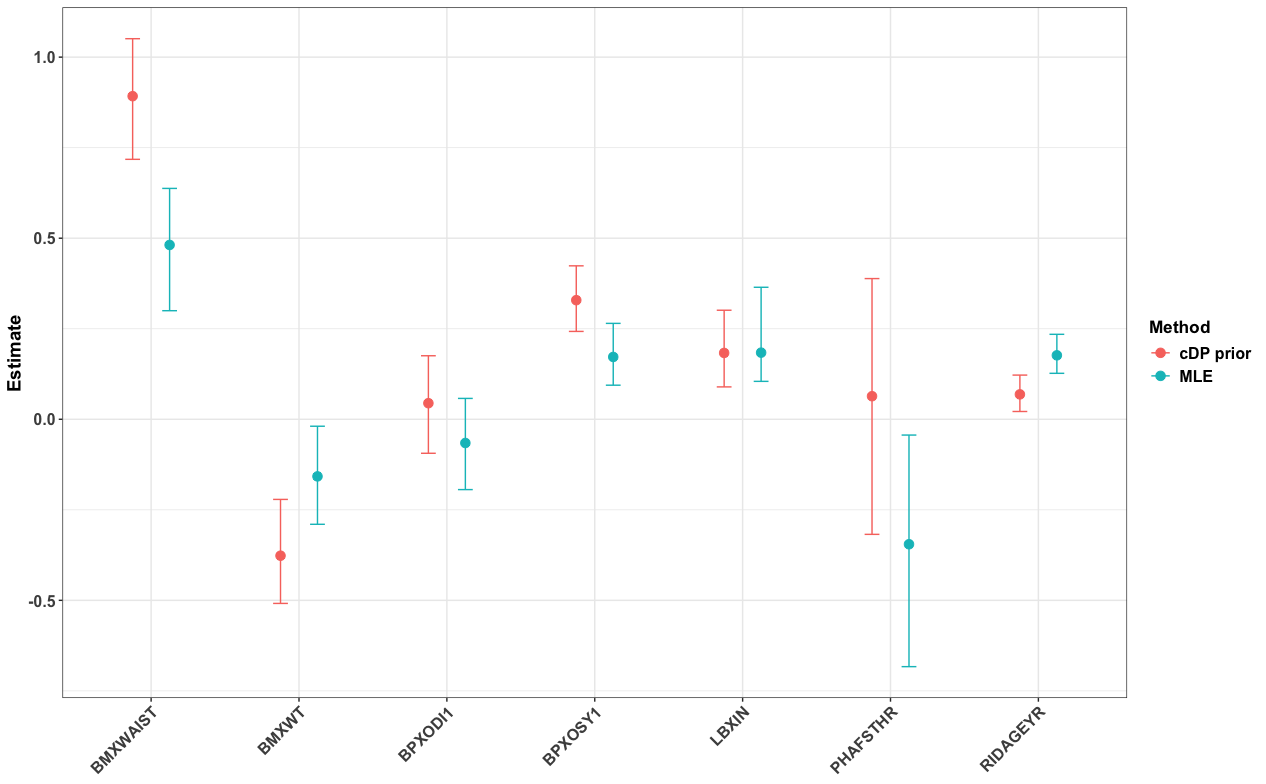}
    \caption{Posterior estimates and 95\% credible intervals for each explanatory variable.}
    \label{fig:post_coef}
\end{figure}

Table~\ref{tab:reg_post_summary} summarizes the estimated values of the coefficients for the regression model with random covariates, along with their associated confidence intervals. These results are visualized in Figure~\ref{fig:post_coef}. For comparison, we report both frequentist estimates, which are obtained via maximum likelihood estimation (MLE) with bootstrap confidence intervals, and posterior estimates based on 1,000 posterior draws under the proposed model.
\begin{table}[H]
\centering
\caption{Posterior summary for regression coefficients}
\begin{tabular}{ll|rrr}
\hline\hline
Explanatory Variable & Method & Estimate & 95\% CI Lower & 95\% CI Upper \\
\hline
BMXWT     & MLE        & $-0.1577$ & $-0.2900$ & $-0.0192$ \\
          & cDP prior  & $-0.3767$ & $-0.5085$ & $-0.2213$ \\
BMXWAIST  & MLE        & $0.4814$  & $0.2998$  & $0.6375$  \\
          & cDP prior  & $0.8922$  & $0.7178$  & $1.0509$  \\
BPXOSY1   & MLE        & $0.1721$  & $0.0941$  & $0.2648$  \\
          & cDP prior  & $0.3289$  & $0.2426$  & $0.4237$  \\
BPXODI1   & MLE        & $-0.0655$ & $-0.1942$ & $0.0576$  \\
          & cDP prior  & $0.0444$  & $-0.0939$ & $0.1754$  \\
LBXIN     & MLE        & $0.1840$  & $0.1047$  & $0.3645$  \\
          & cDP prior  & $0.1832$  & $0.0894$  & $0.3011$  \\
RIDAGEYR  & MLE        & $0.1766$  & $0.1269$  & $0.2346$  \\
          & cDP prior  & $0.0688$  & $0.0215$  & $0.1220$  \\
PHAFSTHR  & MLE        & $-0.3451$ & $-0.6831$ & $-0.0436$ \\
          & cDP prior  & $0.0637$  & $-0.3178$ & $0.3885$  \\
\hline\hline
\end{tabular}
\label{tab:reg_post_summary}
\end{table}

From Table~\ref{tab:reg_post_summary} and Figure~\ref{fig:post_coef}, we can see that waist circumference (\texttt{BMXWAIST}) shows a strong positive association with fasting glucose under both MLE and cDP prior. In contrast, body weight (\texttt{BMXWT}) shows a negative effect in both methods. Diastolic blood pressure (\texttt{BPXODI1}) has wide intervals covering zero, indicating no significant effect. Fasting insulin (\texttt{LBXIN}) shows consistent positive effects under both methods, supporting biological expectations. Hours since last meal (\texttt{PHAFSTHR}) has a negative effect under MLE, but the cDP estimate is near zero, possibly because the scatter plots show no clear trend between fasting time and glucose levels.


\section{ Proofs of Lemma \ref{lm1} and Theorem \ref{th1} }

\begin{proof}[Proof of Lemma \ref{lm1}] Note 
\begin{align*}
Q(d\bfx_n, dF) & = \prod_{i=1}^n F(dx_i) \calD_\alpha(dF) \\
& =  \prod_{i=1}^n F(dx_i) \calD_{\alpha, \xi}(dF) \mu(d\xi) \\
& = \prod_{i=1}^n F(dx_i) \calD_{\alpha} (dF|g(F) = \xi) h(\xi: g, \alpha) \mu(d\xi). 
\end{align*}
Thus,  the measure $Q$ has $(g, \mu)$-disintegration $ (\prod_{i=1}^n F(dx_i) \calD_{\alpha} (dF|g(F) = \xi) h(\xi: g, \alpha),\xi \in \bbR)$. 

On the other hand,  
\begin{align*}
Q(d\bfx_n, dF) & = \calD_{\alpha + \sum_{i=1}^n \delta_{x_i} } (dF) Polya_\alpha(d\bfx_n) \\
& = \calD_{\alpha + \sum_{i=1}^n \delta_{x_i}, \xi} (dF) \mu(d\xi) Polya_\alpha(d\bfx_n) \\
& =  \calD_{\alpha+nF_n} (dF|g(F) = \xi)  h(\xi: g, \alpha+nF_n) Polya_\alpha(d\bfx_n)\mu(d\xi). 
\end{align*}
Thus,  $Q$ also has the $(g, \mu)$-disintegration  $( \calD_{\alpha+nF_n} (dF|g(F) = \xi)  h(\xi: g, \alpha+nF_n) Polya_\alpha(d\bfx_n), \xi \in \bbR)$.  By the uniqueness of the disintegration,  
\begin{align*}
    \prod_{i=1}^n F(dx_i) \calD_{\alpha} (dF|g(F) = \xi) h(\xi: g, \alpha) &= \calD_{\alpha+nF_n} (dF|g(F) = \xi)  h(\xi: g, \alpha+nF_n) \\
    &\phantom{00} \times Polya_\alpha(d\bfx_n), ~ \mu-a.a ~ \xi.
\end{align*}
This completes the proof. 
\end{proof} 

\begin{proof} [Proof of Theorem \ref{th1}]   We calculate $P(dF,  d\xi, d\bepsilon_n )$ first.  
\begin{align*} 
& ~ P(dF,  d\xi, d\bepsilon_n ) \nonumber \\
& = \prod_{i=1}^n F(d\epsilon_i )\pi(dF, d\xi) \nonumber\\
& = \prod_{i=1}^n F(d\epsilon_i ) \calD_{\alpha}(dF| g(F, \xi) = 0) \pi(\theta) d\theta \delta_0(d\nu) \\
& = \calD_{\alpha + n F_{n}^\epsilon} (dF| g(F, \xi) = 0) Polya_\alpha(d\bepsilon_n) \frac{h(0: g(F, \xi), \alpha + n F_{n}^\epsilon)} { h(0: g(F, \xi), \alpha)}  \pi(\theta)d\theta \delta_0(d\nu) \\
& = \calD_{\alpha + n F_{n}^\epsilon} (dF| g(F, \xi) = 0) \Big[ \sum_{k=1}^n \prod_{i=1}^k a(\epsilon_i^*) d\epsilon_i^* p_A(\Pi(\bepsilon_n)) \Big] \frac{h(0: g(F, \xi), \alpha + n F_{n}^\epsilon)} { h(0: g(F, \xi), \alpha)}  \pi(\theta)d\theta \delta_0(d\nu), 
\end{align*}
where $F_n^\epsilon = \frac{1}{n} \sum_{i=1}^n \delta_{\epsilon_i}$. 
The third equality holds due to Lemma \ref{lm1}. 

By the change of variable formula,  we obtain the following.  
\begin{align} \label{eq: joint-F-xi-epsilon} 
& ~ P(dF,  d\xi, d\bfy_n ) \nonumber \\
& = \calD_{\alpha + n F_{n, \xi}} (dF| g(F, \xi) = 0 )\Big[ \sum_{k=1}^n \prod_{i=1}^k a(t(y_i^*, \xi, x_i^*)) | t'(y_i^*, \xi, x_i^*)| dy_i^* p_A(\Pi(\bfy_n))  \Big] \nonumber \\
&\phantom{00} \times \frac{h(0: g(F, \xi),  \alpha + n F_{n, \xi})} { h(0: g(F, \xi), \alpha)}  \pi(\theta)d\theta \delta_0(d\nu) \nonumber \\
& = \calD_{\alpha + n F_{n, \xi}} (dF| g(F, \xi) = 0 )\Big[ \sum_{k=1}^n  \prod_{i=1}^k  b(y_i^*: \xi, x_i^*) dy_i^* p_A(\Pi(\bfy_n))  \Big] \nonumber \\
&\phantom{00} \times \frac{h(0: g(F, \xi),  \alpha + n F_{n, \xi})} { h(0: g(F, \xi), \alpha)}  \pi(\theta)d\theta \delta_0(d\nu).
\end{align}

From equation \eqref{eq: joint-F-xi-epsilon},  we can obtain the posterior,   
\begin{align} \label{eq:post-F-xi} 
\pi(dF, d\xi | \bfy_n ) & \propto \calD_{\alpha + n F_{n, \xi}} (dF| g(F, \xi) = 0 )  \prod_{i=1}^k  b(y_i^*: \xi, x_i^*)  \frac{h(0: g(F, \xi),  \alpha + n F_{n, \xi})} { h(0: g(F, \xi), \alpha)}  \pi(\theta)d\theta \delta_0(d\nu). 
\end{align}
If we integrate out $\nu$ from \eqref{eq:post-F-xi}, we get 
\begin{align} \label{eq:post-F-theta} 
\pi(dF, d\theta | \bfy_n ) & \propto \calD_{\alpha + n F_{n, (\theta, 0)}} (dF| g(F, (\theta, 0)) = 0 ) \nonumber\\
&\phantom{00} \times \prod_{i=1}^k  b(y_i^*: (\theta, 0), x_i^*)  \frac{h(0: g(F, (\theta, 0)),  \alpha + n F_{n, (\theta, 0)})}{ h(0: g(F, (\theta, 0)), \alpha)}  \pi(\theta)d\theta. 
\end{align}

We integrate out  $F$ from \eqref{eq:post-F-xi}, we obtain 
\begin{align}\label{eq:post-theta} 
\pi(d\theta | \bfy_n ) & \propto    \prod_{i=1}^k  b(y_i^*: (\theta, 0), x_i^*)  \frac{h(0: g(F, (\theta, 0)),  \alpha + n F_{n, (\theta, 0)})} { h(0: g(F, (\theta, 0)), \alpha)}  \pi(\theta)d\theta. 
\end{align}
This completes the proof. 
\end{proof} 

\section{ Proofs of Theorems \ref{thm:asymptoticA} and \ref{thm:asymptotic}}

\begin{proof}[Proof of Theorem \ref{thm:asymptoticA}]

    By the assumptions on the convergences of $h(0: g(F, (\theta, 0)), \alpha + n F_{n, (\theta, 0)} )$ and $h_1(\theta, \alpha)$ and the assumption on the bounded condition of $\frac{h(0:g(F, (\theta, 0)),\alpha+ nF_{n, (\theta, 0)})}{ h_1(\theta,\alpha)} \pi(\theta)$, the dominated convergence theorem gives 
    \bea z_A =\int \frac{h(0:g(F, (\theta, 0)),\alpha+ nF_{n, (\theta, 0)})}{ h_1(\theta,\alpha)} \pi(\theta)d\theta \to  \int \frac{h(0: g(F, (\theta, 0)),  n F_{n, (\theta, 0)} )}{ \bar{h}(\theta)} \pi(\theta)d\theta = z_0,
    \eea 
    which also guarantees \eqref{eq:convergenceA}.
    
\end{proof}

\begin{proof}[Proof of Theorem \ref{thm:asymptotic}]

The structure of the proof closely follows that of Theorem 4 in \citeasnoun{miller2021asymptotic}. However, we provide a detailed argument here since our assumptions are slightly weaker than those imposed in their setting.

We show
\bea 
\int_{\bbR^p} |q_n(x) - q_0(x)| dx \to 0,
\eea 
where $q_0(x) = (2\pi)^{-p/2} |H_0|^{1/2} \exp(-x^T H_0^{-1} x/2) $ that represents the density function of $N_p(0,H_0)$.

Define 
\bea 
g_n(x) &=& q_n(x) \exp(-f_n(\theta_n))n^{p/2} z_n\\
&=&  \exp(f_n(\theta_n+x/\sqrt{n}) -f_n(\theta_n)) \pi(\theta_n+x/\sqrt{n}),\\
g_0(x) &=& q_0(x) (2\pi)^{p/2} |H_0|^{-1/2} \pi(\theta_0) \\
&=&
\exp(-\frac{1}{2} x^T H_0 x) \pi(\theta_0),
\eea
and define $a_n$ and $a_0$ to satisfy $a_n g_n(x) = q_n(x)$ and $a_0 g_0(x) = q_0(x)$ as below:
\bea 
a_n &=& \frac{\exp(f_n(\theta_n))}{n^{p/2}z_n}\\
a_0 &=& |H_0|^{1/2} /(2\pi)^{p/2}.
\eea 
Since 
\bea 
\int |q_n(x) - q_0(x)| dx &=&
    \int |a_n g_n(x) - a_0 g_0(x)| dx\\
    &\le& \int |a_n g_n(x) - a_n g_0(x)| dx +\int |a_n g_0(x) - a_0 g_0(x)| dx\\
    &=& |a_n| \int |g_n(x) - g_0 (x)| dx + |a_n - a| \int |g_0(x)| dx,
\eea 
$a_n = ( \int g_n(x) dx)^{-1}$ and $a_0 = ( \int g_0(x) dx)^{-1}<\infty$, 
it suffices to show $\int |g_n(x) - g_0(x)| dx \to 0$, which is to be shown by the generalized dominated convergence theorem. 

First, we show that $g_n$ converges to $g$ pointwisely.
We have 
\bea 
\log (g_n(x)) &=& f_n(\theta_n+x/\sqrt{n}) - f_n(\theta_n) + \log \pi(\theta_n+x/\sqrt{n}) \\
&=& -\frac{1}{2} x^T H_n x + r_n(\theta_n +x/\sqrt{n}) -r_n(\theta_n ) + \log \pi(\theta_n+x/\sqrt{n}) \\
&\to& -\frac{1}{2} x^T H_0 x + \log \pi(\theta_0),~ \textup{as } n\to \infty,
\eea 
since $\pi$ is continuous at $\theta_0$ and $ |r_n(\theta_n +x/\sqrt{n}) -r_n(\theta_n )|\to 0$ by the assumption.
Thus, $g_n(x)\to g(x)$ for each $x$.

We define $G_n(x)$ and $G_0(x)$, considering $|g_n| \leq G_n$ and $G_n(x)\to G_0(x)$, as    
\bea 
G_n(x) &=& \begin{cases}
        \exp\left(-\frac{1}{2}(x- z_n)^T (H_n - c_0 \epsilon I)  (x- z_n)\right)\pi(\theta_0) 2 \exp(M)  & \text{if } |x| < \varepsilon\sqrt{n}, \\
        \exp\left(-n\delta/2\right)\pi\left(\theta_n + \frac{x}{\sqrt{n}}\right)  & \text{if } |x| \geq \varepsilon\sqrt{n},
    \end{cases} \\
    G_0(x) &=& \exp\left(-\frac{1}{2}x^T (H_0 - c_0 \epsilon I) x\right)\pi(\theta_0) 2\exp(M) ,
\eea 
where 
$z_n =  (H_n - c_0\epsilon I)^{-1}c_n $ and 
$\delta = \liminf_n \inf_{\theta \in B_\varepsilon(\theta_n)^c} [ \frac{f_n(\theta_n) -f_n(\theta) }{n}] >0$. Here $\epsilon>0$ is sufficiently small constant, the condition of which is given in this proof. 
Fixing $x \in \mathbb{R}^p$, $|x| < \epsilon\sqrt{n}$ for all $n$ sufficiently large. Since $||z_n||\to 0$ and $H_n\to H_0$ by Condition (C1), $G_n$ converges to $G_0$ pointwisely.

Next, we show $|g_n| \leq G_n$ for all sufficiently large $n$.
Suppose $|x| \geq \varepsilon \sqrt{n}$.
Since $\theta_n + x/\sqrt{n} \in B_\epsilon (\theta_n)^c$ and 
$    \inf_{\theta \in B_\varepsilon(\theta_n)^c} \left(\frac{f_n(\theta_n) - f_n(\theta)}{n}\right) > \delta/2$ for all sufficiently large $n$, 
we obtain 
$ f_n(\theta_n) -f_n\left(\theta_n + \frac{x}{\sqrt{n}}\right) > n\delta/2$ and
\bea 
    g_n(x) &=&\exp(f_n(\theta_n+x/\sqrt{n}) -f_n(\theta_n)) \pi(\theta_n+x/\sqrt{n}) \\
    &\le& \exp\left(-n\delta/2\right)\pi\left(\theta_n + \frac{x}{\sqrt{n}}\right) =G_n(x), ~ |x| \geq \varepsilon \sqrt{n},
\eea 
for all sufficiently large $n$.

Next, suppose $|x| < \varepsilon\sqrt{n}$. 
By setting $\epsilon>0$ small enough, we have $\pi\left(\theta_n + \frac{x}{\sqrt{n}}\right) \leq \sqrt{2}\pi(\theta_0)$ for all sufficiently large $n$. Then, 
\bea 
\exp(f_n(\theta_n+x/\sqrt{n}) -f_n(\theta_n)) &=& \exp(-\frac{1}{2}x^T H_n x  + r_n(\theta_n+x/\sqrt{n}) -r_n(\theta_n) ) \\
&\le& \exp(-\frac{1}{2}x^T H_n x ) \exp(c_0 ||x||^3/(2\sqrt{n})) \exp( x^T c_n +  d_n )  \\
&\le& \exp(-\frac{1}{2}x^T H_n x ) \exp(c_0 \epsilon x^T x /2 +x^T c_n) \exp( d_n )  \\
&\le& \sqrt{2}\exp(-\frac{1}{2}(x-  z_n )^T (H_n - c_0\epsilon I) (x-  z_n ) )  \exp( M )  ,
\eea 
for all sufficiently large $n$, where
the last inequality is satisfied since $\limsup_n d_n\to M$ is satisfied by the assumption.
Thus, 
\bea 
    g_n(x) &=&\exp(f_n(\theta_n+x/\sqrt{n}) -f_n(\theta_n)) \pi(\theta_n+x/\sqrt{n}) \\
    &\le& \exp(-\frac{1}{2}(x-z_n)^T (H_n -c_0\epsilon I ) (x- z_n) )\pi(\theta_0) 2\exp(M) =G_n(x), ~ |x| < \varepsilon \sqrt{n}.
\eea

Finally, we show $\int G_n(x) dx \to \int G_0(x) dx$.
 We have
\bea 
    \int G_n(x) dx &=& \pi(\theta_0)2\exp(M)\int_{B_{\varepsilon \sqrt{n}}(0)} \exp\left(-\frac{1}{2} (x-z_n)^T (H_n - c_0 \epsilon I) (x-z_n)\right) dx \\
    &&+  \int_{B_{\varepsilon \sqrt{n}}(0)^c} \exp(-n\delta/2)\pi\left(\theta_n + \frac{x}{\sqrt{n}}\right)dx.
\eea 
We have 
\bea 
 \int_{B_{\varepsilon \sqrt{n}}(0)^c} \exp(-n\delta/2)\pi\left(\theta_n + \frac{x}{\sqrt{n}}\right)dx&\le&  \exp(-n\delta/2)n^{p/2}  \int_{\bbR^p} \pi\left(\theta_n + t\right)dt\\
 &=&\exp(-n\delta/2)n^{p/2}  ,
\eea 
which converges to $0$.
Next, using the dominated convergence theorem, we show
\bea \int_{B_{\varepsilon \sqrt{n}}(0)} \exp\left(-\frac{1}{2} (x-z_n)^T (H_n - c_0 \epsilon I) (x-z_n)\right) dx \to \int_{\bbR^p} \exp\left(-\frac{1}{2}x^T (H_0 - c_0 \epsilon I) x\right)  dx.
\eea 
For fixed $x$, $\exp\left(-\frac{1}{2} (x-z_n)^T (H_n - c_0 \epsilon I) (x-z_n)\right)  I( x\in B_{\varepsilon \sqrt{n}}(0))$ converges to $\exp\left(-\frac{1}{2}x^T (H_0 - c_0 \epsilon I) x\right) $.
Since for all sufficiently large $n$
\bea 
\exp\left(-\frac{1}{2} (x-z_n)^T (H_n - c_0 \epsilon I) (x-z_n)\right)  I( x\in B_{\varepsilon \sqrt{n}}(0)) \le \exp(-\frac{1}{4} x^T H_0 x),
\eea 
which is integrable. Thus, the dominated convergence theorem gives $\int G_n (x) dx \to \int G_0(x) dx$.



\end{proof}

\section{Proofs of Theorem \ref{thm:quantileAasymp} and Corollary \ref{corr:quantileAasymp}}\label{sec:proofquantileasymp}

\begin{proof}[Proof of Theorem \ref{thm:quantileAasymp}]
    We have
    \bea 
    h(\Delta p: g(F,\theta),\alpha+ nF_n) &=& 
 \frac{\Gamma(\alpha(\mathbb{R}) + n)}{\prod_{j=1}^{k+1} \Gamma(\alpha(S_j(\theta)) + n F_n(S_j(\theta))) } \prod_{j=1}^{k+1} \Delta p_j^{\alpha(S_j(\theta))+ n F_n(S_j(\theta)) -1}      
    \\
    &\to &
 \frac{\Gamma(n)}{\prod_{j=1}^{k+1}\Gamma( n F_n(S_j(\theta))) } \prod_{j=1}^{k+1} \Delta p_j^{n F_n(S_j(\theta))-1}  ,
    \eea 
    as $A\to 0$.

    We decompose $h(\Delta p: g(F,\theta),\alpha)$ as
    \bea 
    h(\Delta p: g(F,\theta),\alpha) &=& \frac{\Gamma(\alpha(\mathbb{R}))}{\prod_{j=1}^{k+1} \Gamma(\alpha(S_j(\theta)))} \prod_{j=1}^{k+1} \Delta p_j^{\alpha(S_j(\theta)) -1}  \\
    &=&h_1(\theta,\alpha) h_2(A),\\
    h_1(\theta,\alpha) &=& \frac{1}{A^{k+1}}\frac{1}{\prod_{j=1}^{k+1}\Gamma(\alpha(S_j(\theta)))} \prod_{j=1}^{k+1} \Delta p_j^{\alpha(S_j(\theta))-1},\\
    h_2(A) &=&A^{k+1}\Gamma(A).
    \eea 
Since \bea 
\frac{1}{A\Gamma(\alpha(S_j(\theta)))} = \frac{\alpha(S_j(\theta))}{A} \frac{1}{\Gamma(\alpha(S_j(\theta)) +1)}  \to \bar\alpha(S_j(\theta))
\eea as $A\to 0$, we have 
    \bea 
h_1(\theta,\alpha)  \to \bar{h}(\theta) = \prod_{j=1}^{k+1} \bar\alpha(S_j(\theta)) \Delta p_j^{-1},~ \textup{as } A\to 0.
    \eea 

    When $F_n(S_j(\theta))>0$ and $A$ is sufficiently small, 
    \bea 
    \frac{A\Gamma(\alpha(S_j(\theta)))}{\Gamma(nF_n(S_j(\theta)) + \alpha(S_j(\theta)))} &\le&\frac{A\Gamma(\alpha(S_j(\theta)))}{\Gamma(nF_n(S_j(\theta)))} \\
    &=& \Gamma(\alpha(S_j(\theta))+1) \frac{A}{\alpha(S_j(\theta))}  \frac{1}{\Gamma(nF_n(S_j(\theta)))} \\
    &\le& 2 \frac{A}{\alpha(S_j(\theta))}  \frac{1}{\Gamma(nF_n(S_j(\theta)))}\\
    &=& 2 \frac{1}{\bar{\alpha}(S_j(\theta))}  \frac{1}{\Gamma(nF_n(S_j(\theta)))}.
    \eea 
    When $F_n(S_j(\theta))=0$, $\frac{A\Gamma(\alpha(S_j(\theta)))}{\Gamma(nF_n(S_j(\theta)) + \alpha(S_j(\theta)))} = A$.
    Thus, when $A$ is sufficiently small, 
\bean 
&&\frac{h(0: g(F, (\theta, 0)), \alpha + n F_{n, (\theta, 0)} )}{ h_1(\theta, \alpha  ) }  \pi(\theta) \nonumber\\
&=& A^{k+1}\prod_{j=1}^{k+1} \frac{\Gamma(\alpha(S_j(\theta)))}{\Gamma(nF_n(S_j(\theta)) + \alpha(S_j(\theta)))} 
\Gamma(\alpha(\bbR)+n)
\prod_{j=1}^{k+1} \Delta p_j^{nF_n(S_j(\theta))} \pi(\theta)  \nonumber\\
&\le& 
  \frac{2}{\prod_{j\in \bm{J}}\bar{\alpha}(S_j(\theta))}\frac{\Gamma(n+1)}{\prod_{j\in \bm{J}}\Gamma(nF_n(S_j(\theta))) } \prod_{j=1}^{k+1} \Delta p_j^{nF_n(S_j(\theta))} \pi(\theta),\label{eq:upper1}
    \eean 
    where $\bm{J} = \{ j\in [k+1] : F_n(S_j(\theta))>0 \} $.
    Thus, since we assume that \eqref{eq:upper1} is integrable, Theorem \ref{thm:asymptoticA} gives \bea  
    \pi(\theta\mid   \bfy_n )  \to
  \frac{1}{z_0\,\prod_{j\in\bm{J}}\bar\alpha(S_j(\theta))
          \Gamma\!\bigl(nF_n(S_j(\theta))\bigr)}
  \prod_{j=1}^{k+1}\Delta p_j^{\,nF_n(S_j(\theta))}\pi(\theta),~ \forall \theta.
    \eea 
    
\end{proof}

\begin{proof}[Proof of Corollary \ref{corr:quantileAasymp}]
    By setting $\bar{\alpha}(-\infty,\theta] = F_n(-\infty,\theta]$, \eqref{eq:upper1} becomes
    \bea 
     2n^{k+1} \frac{n!}{\prod_{j=1}^{k+1} [nF_n(S_j(\theta))]!} \prod_{j=1}^{k+1} \Delta p_j^{nF_n(S_j(\theta))}\pi(\theta),
    \eea 
    then 
    \bea 
     \eqref{eq:upper1}&=&2n^{k+1}\int \frac{n!}{\prod_{j=1}^{k+1} [nF_n(S_j(\theta))]!} \prod_{j=1}^{k+1} \Delta p_j^{nF_n(S_j(\theta))} \pi(\theta) d\theta \\
     &=& 2n^{k+1}\int p_{MN}(nF_n(S_j(\theta)); \Delta p) \pi(\theta) d\theta \\
     &\le& 2n^{k+1} (\sum_{x} p_{MN}(x; \Delta p))^{k+1} \\
     &=& 2n^{k+1},
    \eea 
    where $p_{MN} $ represents the density function of the multinomial distribution.
    Thus, the choice of $\alpha$ makes Condition (QA1) hold, and the limiting distribution \eqref{eq:limitingAquantiles} is proportional to  
    \bea 
     \frac{1}{\prod_{j=1}^{k+1}(nF_n(S_j(\theta))\bigr)!}
  \prod_{j=1}^{k+1}(\Delta p_j)^{nF_n(S_j(\theta))}
  \pi(\theta).
    \eea 
\end{proof}

\section{Proof of Theorem \ref{thm:asympquantiles}}\label{sec:proofquantileasymp2}

In this section, for a distribution function $F$, we denote
\bea 
F(t) &=& F(-\infty,t],~ t\in \bbR\\
F(\theta) &=& (F(\theta_1),\ldots, F(\theta_k))^T,~\theta= (\theta_1,\ldots,\theta_k)^T\in \bbR^k,\\
\Delta F(\theta) &=&  (F(\theta_1),F(\theta_2)-F(\theta_1)\ldots, F(\theta_k)-F(\theta_{k-1}))^T.
\eea

For the proof of Theorem \ref{thm:asympquantiles}, we provide Propositions \ref{prop:LocalMNPmfExpansion} and \ref{prop:condition2}.
\begin{proposition}\label{prop:LocalMNPmfExpansion}
Suppose the same setting of Theorem \ref{thm:asympquantiles}.
Then, there exist sequences \( C_n \in \mathbb{R} \) and \( B_n \in \mathbb{R}^{m \times m} \), which do not depend on $\theta$, such that
\[
\log \frac{    h(\Delta p: g(F,\theta),\alpha+ nF_n)}{    h(\Delta p: g(F,\theta),\alpha)}=
-\,\frac{n}{2} \, (\theta - q_n)^\top B_n (\theta - q_n)
+ C_n
+ R(\theta),
\]
and $B_n \to  T^T H T$,
where $T\in\mathbb{R}^{m\times m}$ is the lower bidiagonal matrix whose nonzero entries are \(T_{j,j} = f(\theta_{0,j})\), \(T_{j,j-1} = -f(\theta_{0,j-1})\),
$H = \textup{diag}( \frac{1}{[\Delta p]_1} ,\ldots,  \frac{1}{[\Delta p]_{k}}  ) + (\frac{1}{1-p_k}) \bm{1}\bm{1}^T$, and the remainder term \( R(\theta) \) satisfies:
\begin{enumerate}
  \item For any fixed \( x \in \mathbb{R}^m \),  
  \[
  R\big(q_n + x/\sqrt{n}\big) - R(q_n) \to 0,
  \quad \text{as } n \to \infty,
  \]
  \item For all \( x \in B_{\epsilon \sqrt{n}}(0) \), the following bound holds:
  \[
  \big| R(q_n + x/\sqrt{n}) - R(q_n) \big|
  \le
  \frac{\|x\|^3}{\sqrt{n}} + x^\top c_n + d_n,
  \]
  where \( c_n \in \mathbb{R}^m \) and \( d_n \in \mathbb{R} \) are sequences satisfying \( ||c_n|| \to 0 \), \( \limsup_{n\to \infty}d_n <\infty \).
\end{enumerate} 
\end{proposition}
The proof of Proposition \ref{prop:LocalMNPmfExpansion} is given in Section \ref{ssec:prop1}.

\begin{proposition}\label{prop:condition2}
    Suppose the same setting of Theorem \ref{thm:asympquantiles}.
    For any $\epsilon$, there exists a constant $C>0$ such that 
    \bea 
    P\Big[\inf_{\theta : |\theta- q_n| >\epsilon} \frac{f_n(q_n) -f_n(\theta)}{n}  > C \Big] \to 1,
    \eea 
    where $f_n(\theta) = \frac{    h(\Delta p: g(F,\theta),\alpha+ nF_n)}{    h(\Delta p: g(F,\theta),\alpha)}$.
\end{proposition}
The proof of Proposition \ref{prop:condition2} is given in Section \ref{ssec:prop2}.

\begin{proof}[Proof of Theorem \ref{thm:asympquantiles}]
    
    The conditions (C1) and (C2) are satisfied by Propositions \ref{prop:LocalMNPmfExpansion} and \ref{prop:condition2}, respectively, which completes the proof of Theorem \ref{thm:asympquantiles}.
\end{proof}

\subsection{Proof of Proposition \ref{prop:LocalMNPmfExpansion}}\label{ssec:prop1}

To prove Proposition \ref{prop:LocalMNPmfExpansion}, we provide Lemmas \ref{lemma:multinomial}-\ref{lem:logl}.

\begin{lemma}\label{lemma:multinomial}
Let \( \mathbf{X} = (X_1, \dots, X_{k}) \sim \mathrm{Multinomial}(n, (p_1, \dots, p_k) ) \), where \( p_i > 0 \), \( \sum_{i=1}^{k} p_i = 1 \), and \( X_k = n - \sum_{i=1}^{k-1} X_i \), \( p_k = 1 - \sum_{i=1}^{k-1} p_i \).
Let $z_i = \frac{X_i - n p_i}{\sqrt{n p_i}}$. 
\bea 
\log p(\bm{X}) &=&  -\frac{k-1}{2}\log (2\pi) -\sum_{i=1}^k \frac{1}{2}\log p_i  \\
&&-\frac{n}{2} ( X_{-k}/n - p_{-k} + \tilde{H}^{-1}\tilde{b} )^T \tilde{H} ( X_{-k}/n - p_{-k}+ \tilde{H}^{-1}\tilde{b} ) + \frac{n}{2}\tilde{b}^T \tilde{H}^{-1} \tilde{b} + R\\
|R|&\le & \frac{1}{12n}+\sum_{i=1}^k\Big[\frac{1}{12X_i} -\frac{ ( \sqrt{n}(X_i/n - p_i)  )^3}{3 (1+x_i^*)^3 \sqrt{n}p_i} (1 +  \frac{\sqrt{n}(X_i/n - p_i)}{ \sqrt{n}p_i} + 1/(2np_i) ) \Big],
\eea
where $\tilde{H} = \textup{diag}( \frac{1}{p_1} - \frac{1}{2np_1^2},\ldots,  \frac{1}{p_{k-1}} - \frac{1}{2np_{k-1}^2} ) + (\frac{1}{p_k}-\frac{1}{2np_k^2}) \bm{1}\bm{1}^T$, $\tilde{b} = \frac{1}{2}( \frac{1}{n\sqrt{p_1}} - \frac{\sqrt{p_1}}{n^{3/2} p_k} , \ldots, \frac{1}{n\sqrt{p_{k-1}}} - \frac{\sqrt{p_{k-1}}}{n^{3/2} p_k})^T$ and $|x_i^*| \le |\frac{X_i}{np}-1| $.
\end{lemma}

\begin{proof}

The log probability of the multinomial distribution is 
\[
\log p(\mathbf{X}) = \log \frac{n!}{\prod_{i=1}^{k} X_i!} + \sum_{i=1}^{k} X_i \log p_i.
\]
Using Stirling’s approximation
\bean\label{eq:Stirling} 
\log(m!) = m \log m - m + \frac{1}{2} \log(2\pi m)  +
R(m), \quad | R(m)| < \frac{1}{m}, 
\eean 
we obtain
\bea 
\log p(\mathbf{X}) 
&=& n \log n - n - \sum_{i=1}^{k} (X_i \log X_i - X_i) + \sum_{i=1}^{k} X_i \log p_i \\
&& + \frac{1}{2} \log(2\pi n) - \frac{1}{2} \sum_{i=1}^{k} \log(2\pi X_i) + R(n) - \sum_{i=1}^{k} R(X_i)\\
&=& \Big(n + \frac{1}{2}\Big)\log n  - \frac{k-1}{2} \log (2\pi) + R(n)\\
&&+\sum_{i=1}^k \Big[ -X_i \log X_i  +X_i\log p_i -\frac{1}{2} \log X_i +R(X_i) \Big] .
\eea 
Since $ X_i = n p_i + \sqrt{n p_i} z_i$ and $\log (1+x) = x -x^2/2 +x^3/\{ 3 (1+ x^*)^3 \} $ for some $x^*$ satisfying $|x^*| <|x|$, we have 
\bea 
\log (X_i) &=&
\log (np_i) + \log (1+\frac{z_i}{\sqrt{np_i}})\\
&=& \log(n p_i) + \frac{z_i}{\sqrt{n p_i}} - \frac{z_i^2}{2 n p_i} + \frac{z_i^3}{3 (1+z_i^*)^3 (np)^{3/2}},
\eea 
where $|z_i^*|<  |z_i|/\sqrt{n p_i}$.
Then,
\[
\begin{aligned}
X_i \log X_i 
&= (n p_i + \sqrt{n p_i} z_i) \left( \log(n p_i) + \frac{z_i}{\sqrt{n p_i}} - \frac{z_i^2}{2 n p_i} + \frac{z_i^3}{3 (1+z_i^*)^3 (np)^{3/2}} \right) \\
&= n p_i \log(n p_i) + \sqrt{n p_i} z_i \log(n p_i) + \sqrt{n p_i} z_i + z_i^2 - \frac{z_i^2}{2} \\
&\quad + \frac{z_i^3}{3 (1+z_i^*)^3 (np)^{1/2}} (1 + z_i/(np)^{1/2} ).
\end{aligned}
\]
We also have $\sum_{i=1}^{k} X_i \log p_i = n \sum_{i=1}^{k} p_i \log p_i + \sum_{i=1}^{k} \sqrt{n p_i} z_i \log p_i$.
Thus,
\bea 
&&\sum -X_i\log X_i + X_i \log p_i -\frac{1}{2}\log (X_i)\\
&=&\sum_{i=1}^k \Big[- n p_i \log(n p_i) - \sqrt{n p_i} z_i \log(n p_i) - \sqrt{n p_i} z_i - \frac{z_i^2}{2} \\
&&+ np_i \log(p_i) + \sqrt{n p_i} z_i\log p_i  -\frac{1}{2} \log (np_i) -\frac{z_i}{2\sqrt{np_i}} + \frac{z_i^2}{4np_i} \\
&&- \frac{z_i^3}{3 (1+z_i^*)^3 (np_i)^{1/2}} (1 + z_i/(np_i)^{1/2}  + 1/(2np_i) )  \Big]\\
&=&  -(n+1/2)\log n  + \frac{1}{2}\sum_{i=1}^k \Big[ \log p_i  -\frac{z_i}{2\sqrt{np_i}}  - \frac{z_i^2}{2} + \frac{z_i^2}{4np_i} \\
&&-\frac{z_i^3}{3 (1+z_i^*)^3 (np_i)^{1/2}} (1 + z_i/(np_i)^{1/2}  + 1/(2np_i) ) \Big],
\eea 
where the last equality is satisfied by the condition $\sum_{i=1}^{k} \sqrt{n p_i} z_i = 0 $, which gives
\bea 
\log p(\bm{X}) &=& -\frac{k-1}{2}\log (2\pi) -\sum_{i=1}^k \frac{1}{2}\log p_i - \frac{z_i}{2\sqrt{np_i}}  - \frac{z_i^2}{2} + \frac{z_i^2}{4np_i} \\
&&+ R(n)  +\sum_{i=1}^k\Big[ R(X_i) -\frac{z_i^3}{3 (1+z_i^*)^3 (np_i)^{1/2}} (1 + z_i/(np_i)^{1/2}  + 1/(2np_i) ) \Big],
\eea

Next, we express $\log p(\bm{X})$ as the function of $z_{-k} = (z_1,\ldots, z_{k-1})^T$ using
\[
\sum_{i=1}^{k} \sqrt{n p_i} z_i = 0 \quad \Longrightarrow \quad z_k = - \sum_{i=1}^{k-1} \frac{\sqrt{p_i}}{\sqrt{p_k}} z_i.
\]
We have
\bea 
\sum_{i=1}^k z_i/(2\sqrt{np_i}) &=& b^T z_{-k},\\
b &=&  (\frac{1}{2 \sqrt{n p_1}} - \frac{p_1}{2 n p_k},\ldots,\frac{1}{2 \sqrt{n p_{k-1}}} - \frac{p_{k-1}}{2 n p_k}  )^T,
\eea 
and
\bea 
\sum_{i=1}^k \Big[z_i^2 - \frac{z_i^2}{2np_i} \Big]
 &=& \frac{1}{2} z_{-k}^T H z_{-k},\\
 H&=&\textup{diag}(1-\frac{1}{2np_1},\ldots,1-\frac{1}{2np_{k-1}}  ) + (1-\frac{1}{2np_k} ) vv^T  \\
 v &=& (\sqrt{p_1/p_k},\ldots, \sqrt{p_{k-1}/p_k})^T .
\eea
We have
\bea 
z_{-k}^T H z_{-k} &=& n(X_{-k}/n - p_{-k} )^T  \tilde{H}   (X_{-k}/n - p_{-k} ),\\
b^T z_{-k} &=& n \tilde{b}^T  (X_{-k}/n - p_{-k}),
\eea 
where $\tilde{H} = \textup{diag}( \frac{1}{p_1} - \frac{1}{2np_1^2},\ldots,  \frac{1}{p_{k-1}} - \frac{1}{2np_{k-1}^2} ) + (\frac{1}{p_k}-\frac{1}{2np_k^2}) \bm{1}\bm{1}^T$ and $\tilde{b} = \frac{1}{2}( \frac{1}{n\sqrt{p_1}} - \frac{\sqrt{p_1}}{n^{3/2} p_k} , \ldots, \frac{1}{n\sqrt{p_{k-1}}} - \frac{\sqrt{p_{k-1}}}{n^{3/2} p_k})^T$.
Thus, 
\bea 
\log p(\bm{X}) &=& -\frac{k-1}{2}\log (2\pi) -\sum_{i=1}^k \frac{1}{2}\log p_i  - b^T  z_{-k} - \frac{1}{2} z_{-k}^T \tilde{H} z_{-k} + R\\
&=&  -\frac{k-1}{2}\log (2\pi) -\sum_{i=1}^k \frac{1}{2}\log p_i \\ 
&&-\frac{n}{2} ( X_{-k}/n - p_{-k} + \tilde{H}^{-1}\tilde{b} )^T \tilde{H} ( X_{-k}/n - p_{-k}+ \tilde{H}^{-1}\tilde{b} ) + \frac{n}{2}\tilde{b}^T \tilde{H}^{-1} \tilde{b} + R\\
|R|&\le & \frac{1}{12n}+\sum_{i=1}^k\Big[\frac{1}{12X_i} -\frac{z_i^3}{3 (1+z_i^*)^3 (np_i)^{1/2}} (1 + z_i/(np_i)^{1/2}  + 1/(2np_i) ) \Big].
\eea 
By the change of variables, we obtain
\bea 
|R|&\le & \frac{1}{12n}+\sum_{i=1}^k\Big[\frac{1}{12X_i} -\frac{ ( \sqrt{n}(X_i/n - p_i)  )^3}{3 (1+x_i^*)^3 \sqrt{n}p_i} (1 +  \frac{\sqrt{n}(X_i/n - p_i)}{ \sqrt{n}p_i} + 1/(2np_i) ) \Big],
\eea
where $|x_i^*| \le |\frac{X_i}{np}-1| $.

\end{proof}

\begin{lemma}
\label{lem:LocalLinearExpansionLipschitz}
Suppose $X_1,\ldots, X_n\sim F$ independently, where $F$ is a distribution function on $\bbR$. Let $\theta_0\in\mathbb{R}$ be a fixed value and assume 
\begin{enumerate}
\item[(i)] \(F\) is differentiable at \(\theta_0\) with \(f(\theta_0)=F'(\theta_0)>0\),
\item[(ii)] \(f(\theta)\) is Lipschitz continuous in a neighborhood of \(\theta_0\), i.e., there exists \(L>0\) such that
  \[
  \bigl\lvert f(\theta) - f(\theta')\bigr\rvert 
  \;\le\;
  L\,|\theta-\theta'|
  \quad
  \text{for all \(\theta,\theta'\) near \(\theta_0\).}
  \]
\end{enumerate}
Let $p = F(\theta_0)$ and $q_n = F_n^{-1} (p)  = \inf\{ x \in \bbR : F_n(x) \le p \}$, where $F_n(t) 
\;=\; 
\frac{1}{n}\sum_{i=1}^n  I(X_i \le t) $.
If there is no pair $i$ and $j$ such that $X_j = X_i$, then
\bea 
\sqrt{n} \{ F_n(q_n) - F(q_n)\}  = \sqrt{n} f(q_n) (\theta_0 - q_n) + R_n,
\eea 
where $R_n$ satisfies $|R_n|\le 1/\sqrt{n} + L\sqrt{n}|q_n - \theta_0|^2$.

\end{lemma}

\begin{proof}
We have 
\bea 
\sqrt{n} \{ F_n(q_n) - F(q_n)\} 
&=& \sqrt{n} \{ F_n(q_n) - F(\theta_0) \} + \sqrt{n}\{F(\theta_0) - F(q_n)\}
\eea 
Since there is no pair $i$ and $j$ such that $X_j = X_i$ by the assumption, we have
\bea 
|F(\theta_0) -F_n(q_n)| \le   \frac{1}{n}.
\eea 
By a first‐order expansion plus a Lipschitz remainder:
\[
F(q_n)  - F(\theta_0)
\;=\;
f(q_n)\,\bigl( q_n- \theta_0 \bigr)
\;+\;
\bigl[f(\theta^*) - f(q_n)\bigr]\,(q_n-\theta_0),
\]
for some point \(\theta^*\) between \(\theta_0\) and \(q_n\).  
Since \(f\) is \(L\)-Lipschitz near \(\theta_0\), 
\(\bigl\lvert f(\theta^*) - f(\theta_0)\bigr\rvert 
\le L\,\bigl\lvert \theta^*-\theta_0\bigr\rvert.\)
Hence 
\[
 F(\theta_0)-F(q_n)
\;=\;
f(q_n)\,\bigl( \theta_0- q_n\bigr)
\;+\;
O\!\Bigl((\theta_0- q_n)^2\Bigr).
\]
Collecting the equations completes the proof.

\end{proof}

\begin{lemma}\label{lem:shifted-quadratic}
Suppose the same setting as in Theorem~\ref{thm:asympquantiles}. 
Assume there exists a unique vector \( \theta_0 = (\theta_{0,1}, \ldots, \theta_{0,k})^\top \in \mathbb{R}^k \) such that
\[
F(\theta_{0,j}) = p_j, \quad \text{for } j = 1, \ldots, k.
\]
Let \( b \in \mathbb{R}^k \) and \( A \in \mathbb{R}^{k \times k} \) denote a fixed vector and matrix, respectively, that do not depend on \( \theta \).
Then,
\bea 
(\Delta F_n (\theta) - \Delta p_{-k} -b)^T A (\Delta F_n (\theta) - \Delta p_{-k} -b) =
(\theta- q_n)^T T^T A T (\theta- q_n) +  r(\theta) +c_n,
\eea 
where 
$T\in\mathbb{R}^{k\times k}$ be the lower bidiagonal matrix whose nonzero entries are \(T_{j,j} = f(q_{n,j})\), \(T_{j,j-1} = -f(q_{n,j-1})\), $q_n = (q_{n,1},\ldots, q_{n,k})^T$, $q_{n,j} =  F_n^{-1} (\tau_j)  = \inf\{ x \in \bbR : F_n(x) \le \tau_j \} $,   
$c_n$ does not depend on $\theta$, and 
$r(\theta) $ satisfies 
\bea 
|r(\theta)| &\le& C\Big(|| \theta - q_n||^3 + || \theta - q_n|| \Big[
\frac{1}{n} + ||b|| + \frac{|| \Delta G(\theta,\theta_0)|| + || \Delta G(q_n,\theta_0)|| }{\sqrt{n}}  \Big]\\
&&~~~~~~~~~~~~~~~~~~~~~~~~~~~~~~~~+ 
\frac{|| \Delta G(\theta,\theta_0)||^2 + || \Delta G(q_n,\theta_0)||^2 }{n}
\Big),
\eea 
for some positive constant $C$, where $G_n(\theta_j) = \sqrt{n}(F_n(\theta_j ) -F(\theta_j))$ and 
$\Delta G_n(\theta_j,\theta_{0,j}) = G_n(\theta_j) - G_n(\theta_{0,j})$.

\end{lemma}

\begin{proof}
Let $\theta_0 = (\theta_{0,1},\ldots ,\theta_{0,k})^T$.
Since $F(\theta_{0,j})=\tau_j$,
\bea 
&&[\Delta F_n(\theta)-  \Delta \tau]_j \\
&=&[F_n(\theta_j) - F_n(\theta_{j-1})] -  [F(\theta_{0,j}) - F(\theta_{0,j-1})] \\ &=& [F(\theta_j) -F(\theta_{0,j})] - [F(\theta_{j-1}) - F(\theta_{0,j-1})] + G_n(\theta_j)/\sqrt{n} - G_n(\theta_{j-1})/\sqrt{n}\\
&=& [F(\theta_j) -F(\theta_{0,j})] - [F(\theta_{j-1}) - F(\theta_{0,j-1})] +  G_n(\theta_{0,j})/\sqrt{n}
\\ &&- G_n(\theta_{0,j-1})/\sqrt{n} + \Delta G_n(\theta_j ,\theta_{0,j})/\sqrt{n} -\Delta G_n(\theta_{j-1} ,\theta_{0,j-1})/\sqrt{n}
\\
&=& [F(\theta_j) -F(\theta_{0,j}) + \frac{G_n(\theta_{0,j}) + \Delta G_n(\theta_j,\theta_{0,j})}{\sqrt{n}} ]\\
&&- [F(\theta_{j-1}) - F(\theta_{0,j-1}) + \frac{G_n(\theta_{0,j-1}) + \Delta G_n(\theta_{j-1},\theta_{0,j-1})}{\sqrt{n}} ] \\
&=& [F(\theta_j) -F(\theta_{0,j}) + \frac{G_n(q_{n,j})}{\sqrt{n}} + \frac{\Delta G_n(\theta_{0,j},q_{n,j})+ \Delta G_n(\theta_j,\theta_{0,j})}{\sqrt{n}} ]\\
&&- [F(\theta_{j-1}) - F(\theta_{0,j-1}) + \frac{G_n(q_{n,j-1})}{\sqrt{n}} + \frac{\Delta G_n(\theta_{0,j-1},q_{n,j-1})+ \Delta G_n(\theta_{j-1},\theta_{0,j-1})}{\sqrt{n}} ] 
\eea

For each \(j=1,\dots,k\), expand \(F(\theta_j)\) around $q_{n,j}$, 
\bean\label{Fthetaq1} 
F(\theta_j)
\;=\;
F(q_{n,j})
\;+\;
f(q_{n,j})\,\bigl(\theta_j - q_{n,j}\bigr)
\;+\;
R_j(\theta_j),
\eean
where
\[
R_j(\theta_j)
\;=\;
\bigl[f(\theta_j^*) - f(q_{n,j})\bigr]\;\bigl(\theta_j - q_{n,j}\bigr)
\]
for some \(\theta_j^*\) on the line segment between \(\theta_j\) and $q_{n,j}$.  
Since \(f\) is \(L\)-Lipschitz, 
\[
|R_j(\theta_j)|
\;\le\;
L \,\bigl|\theta_j -q_{n,j}\bigr|^2.
\]
Likewise, we have
\bean\label{Fthetaq2}
F(\theta_{0,j}) = F(q_{n,j}) + f(q_{n,j}) ( \theta_{0,j} - q_{n,j}) + R_{0,j},
\eean where $R_{0,j} = \bigl[f(\theta_{0,j}^*) - f(q_{n,j})\bigr]\;\bigl(\theta_{0,j} - q_{n,j}\bigr)$ for some $\theta_{0,j}$ on the line segment between $\theta_{0,j}$ and $q_{n,j}$.
Collecting \eqref{Fthetaq1} and \eqref{Fthetaq2}, we obtain 
\bea 
F(\theta_j) - F(\theta_{0,j}) = f(q_{n,j}) (\theta_j - \theta_{0,j}) + R_j(\theta_j)- R_{0,j}.
\eea 
By Lemma \ref{lem:LocalLinearExpansionLipschitz}, we have
\bea 
G_n(q_{n,j})/\sqrt{n}
&=& f(q_{n,j}) (\theta_{0,j}- q_{n,j}) + R_{0,j}^{(2)} ,
\eea 
where $ |R_{0,j}^{(2)}| \le 1/n + L|q_{n,j} -\theta_{0,j}|^2$.
So, we obtain 
\bea 
[F(\theta_j) - F(\theta_{0,j}) +G_n(q_{n,j})/\sqrt{n}] &=& 
f(q_{n,j}) (\theta_j - q_{n,j}) +
R_j(\theta_j)- R_{0,j} +R_{0,j}^{(2)} .
\eea 
Thus, 
\bea 
\Delta F_n(\theta) - \Delta F(\theta_0)
&=&
T(\theta - q_n + \tilde{R}(\theta) + \Delta \tilde{G}_n(\theta)/\sqrt{n}  + \tilde{R}_0),
\eea 
where $\tilde{R}(\theta) = [R_1(\theta_1)/f(q_{n,1}), \ldots,  R_m(\theta_m)/f(q_{n,k})]$, \\
$ \Delta \tilde{G}_n (\theta) = [\{ \Delta G_n(\theta_{0,1},q_{n,1})+ \Delta G_n(\theta_1,\theta_{0,1})\}/f(q_{n,1}), \ldots, \{ \Delta G_n(\theta_{0,k},q_{n,k})+ \Delta G_n(\theta_k,\theta_{0,k}) \}/f(q_{n,k})]$, and 
$ \tilde{R}_0 = (\tilde{R}_{0,1},\ldots, \tilde{R}_{0,k})^T$ with $\tilde{R}_{0,j} = (R_{0,j}+R_{0,j}^{(2)} )/f(q_{n,j})$.

\bea 
&&(\Delta F_n(\theta) - \Delta F(\theta_0) -b)^T A  (\Delta F_n(\theta) - \Delta F(\theta_0) - b)\\
&=& (\theta- q_n)^T T^T A T (\theta- q_n)  + 2(\theta- q_n)^T T^T A   T \tilde{R}(\theta) \\
&&+ 2(\theta- q_n +\tilde{R}(\theta)  )^T T^T A   (T \tilde{R}_0 + \Delta \tilde{G}_n(\theta)/\sqrt{n} - b)\\
&&+ \tilde{R}(\theta)^T T^T A T \tilde{R}(\theta)  +
  \Delta \tilde{G}_n(\theta)^T T^T A T \Delta \tilde{G}_n(\theta)/n \\
  &&+
(T \tilde{R}_0  - b)^T A   (T \tilde{R}_0 - b)\\
&=& (\theta- q_n)^T T^T A T (\theta- q_n) +  r(\theta) +c_n,
\eea 
where $c_n= (T \tilde{R}_0 + T \tilde{G}_n/\sqrt{n} - b)^T A   (T \tilde{R}_0 + T \tilde{G}_n/\sqrt{n} - b)$ that does not depend on $\theta$.
Since $||\tilde{R}_0 ||_2 \le ||q_n - \theta_0||^2 + 1/n$ and $||\tilde{R}(\theta)||_2 \le ||\theta- q_n||^2$,
$r(\theta)$ satisfies 
\bea 
|r(\theta)| &\le& C\Big(|| \theta - q_n||^3 + || \theta - q_n|| \Big[
\frac{1}{n} + ||b|| + \frac{|| \Delta G(\theta,\theta_0)|| + || \Delta G(q_n,\theta_0)|| }{\sqrt{n}}  \Big]\\
&&~~~~~~~~~~~~~~~~~~~~~~~~~~~~~~~~+ 
\frac{|| \Delta G(\theta,\theta_0)||^2 + || \Delta G(q_n,\theta_0)||^2 }{n}
\Big),
\eea 
for some positive constant $C$.

\end{proof}

\begin{lemma}\label{lem:logl}
Let \( \alpha_i > 0 \) for \( i = 1, \dots, m+1 \), and let \( k_i \geq 0 \) be non-negative integers satisfying \( \sum_{i=1}^{m+1} k_i = n \). 
Define \( \alpha(\mathbb{R}) = \sum_{i=1}^{m+1} \alpha_i \).
\bea 
\sum_{j=1}^{m+1}\log \Big[ \frac{\Gamma(k_j+1)}{\Gamma(\alpha_j+k_j)} \Big] = \sum_{j=1}^{m+1} \Big[ k_j  \log \left( \frac{ k_j + 1 }{ \alpha_j + k_j } \right) 
+ \frac{1}{2}\log (k_j+1)  - (\alpha_j -\frac{1}{2}) \log (\alpha_j + k_j) \Big] + R_n
\eea 
where the remainder satisfies
\[
\left| R_n \right| \leq \frac{1}{12}\sum_{i=1}^{m+1} \Big( \frac{1}{k_i+1} +\frac{1}{\alpha_i + k_i} \Big).
\]
\end{lemma}

\begin{proof}

Applying the log-Gamma Stirling expansion (\citeasnoun{robbins1955remark}):
\[
\log \Gamma(x) = \left( x - \tfrac{1}{2} \right) \log x - x + \tfrac{1}{2} \log(2\pi) + R(x), \quad |R(x)| \leq \frac{1}{12x},
\]
We obtain
\[
\begin{aligned}
\log k_j! &= \left( k_j + \tfrac{1}{2} \right) \log(k_j+1) - (k_j+1) + \tfrac{1}{2} \log(2\pi) + R(k_j+1), \\
\log \Gamma(\alpha_j + k_j) &= \left( \alpha_j + k_j - \tfrac{1}{2} \right) \log (\alpha_j + k_j) - (\alpha_j + k_j) + \tfrac{1}{2} \log(2\pi) + R(\alpha_j + k_j).
\end{aligned}
\]
Since $\sum \big(  - (k_j + 1) + (\alpha_j + k_j) \big)= \alpha(\bbR)-m-1$, 
\bea 
\sum_{j=1}^{m+1}\log \Big[ \frac{\Gamma(k_j+1)}{\Gamma(\alpha_j+k_j)} \Big]&=&  \alpha(\bbR)-m-1 + \sum_{j=1}^{m+1} T_j +R_n,\\
T_j &=& - (\alpha_j + k_j - \tfrac{1}{2}) \log (\alpha_j + k_j) 
+ \left( k_j + \tfrac{1}{2} \right) \log(k_j + 1)\\
&=& k_j  \log \left( \frac{ k_j + 1 }{ \alpha_j + k_j } \right) 
+ \frac{1}{2}\log (k_j+1)  - (\alpha_j -\frac{1}{2}) \log (\alpha_j + k_j)  ,
\eea 
where
\[
|R_n| \le \frac{1}{12}\sum_{i=1}^{m+1} \Big( \frac{1}{k_i+1} +\frac{1}{\alpha_i + k_i} \Big). 
\]

\end{proof}
\begin{lemma}\label{lemma:AsyEquicontProof}
Let \( G_n(t) := \sqrt{n}\,\bigl[F_n(t) - F(t)\bigr] \), where \( F_n \) is the empirical cumulative distribution function based on i.i.d. samples from a continuous distribution \( F \).  
If \( P\bigl( |\theta_n^{(1)} - \theta_n^{(2)}| > \delta \bigr) \to 0 \) as \( n \to \infty \), then
\[
G_n(\theta_n^{(1)}) - G_n(\theta_n^{(2)})
\;\xrightarrow{p}\; 0
\quad
\text{as } n \to \infty.
\]
\end{lemma}

\begin{proof}
Let \( \delta > 0 \) be arbitrary. Then for any \( \epsilon > 0 \), we have
\begin{align*}
P\bigl( |G_n(\theta_n^{(1)}) - G_n(\theta_n^{(2)})| > \epsilon \bigr)
&\le 
P\bigl( |\theta_n^{(1)} - \theta_n^{(2)}| > \delta \bigr)
+ 
P\left( \sup_{|t - s| \le \delta} |G_n(t) - G_n(s)| > \epsilon \right).
\end{align*}
Taking \(\limsup_{n \to \infty}\) on both sides yields
\bean\label{eq:empiricalbound}
\limsup_{n \to \infty} P\bigl( |G_n(\theta_n^{(1)}) - G_n(\theta_n^{(2)})| > \epsilon \bigr)
\le 
\limsup_{n \to \infty} 
P\left( \sup_{|t - s| \le \delta} |G_n(t) - G_n(s)| > \epsilon \right).
\eean 

A classical result in empirical process theory (see, e.g., Pollard (1984), or van der Vaart and Wellner (1996), Chapters 2--3) gives that \( \{G_n(t)\} \) is asymptotically equicontinuous in probability, meaning
\[
\forall \varepsilon > 0,\quad
\lim_{\delta \to 0} \limsup_{n \to \infty}
P\left(
\sup_{|t - s| < \delta} |G_n(t) - G_n(s)| > \varepsilon
\right)
= 0.
\]

Since the right-hand side of \eqref{eq:empiricalbound} can be made arbitrarily small by taking \( \delta \to 0 \), we conclude that the left-hand side must converge to 0 as well. This proves the result.
\end{proof}

\begin{proof}[Proof of Proposition \ref{prop:LocalMNPmfExpansion}]

We have 
\bea 
\exp(f_n(\theta)) &=& p_{MN}(n\Delta F_n(\theta);  \Delta p) l(n\Delta F_n(\theta),\Delta \alpha(\theta))\frac{\Gamma(A+n)}{\Gamma(A)\Gamma(n+1)},\\
l(n\Delta F_n(\theta),\Delta \alpha(\theta)) &=& 
\prod_{j=1}^{k+1}\frac{\Gamma\left( [\Delta \alpha(\theta)]_j \right)  \Gamma( n[\Delta F_n(\theta  )]_j+1)}{ \Gamma\left( [\Delta\alpha(\theta)]_j + n[\Delta F_n(\theta )]_j \right) }.
\eea 
Lemmas \ref{lemma:multinomial} and \ref{lem:shifted-quadratic} give
\bea 
\log p_{MN}(n\Delta F_n(\theta);  \Delta p) &=&- \frac{n}{2}(\theta - q_n)^T T_n^T\tilde{H}_n T_n(\theta - q_n) + C_{n,1} +  nR_{n,1}(\theta) \\
|R_{n,1}(\theta)| &\lesssim &   || \theta - q_n||^3 + || \theta - q_n|| \Big[
\frac{1}{n} + ||\tilde{b}_n|| + \frac{|| \Delta G(\theta,\theta_0)|| + || \Delta G(q_n,\theta_0)|| }{\sqrt{n}}  \Big]\\
&&+
\frac{|| \Delta G(\theta,\theta_0)||^2 + || \Delta G(q_n,\theta_0)||^2 }{n},
\eea
where $T_n\in\mathbb{R}^{k\times k}$ is the lower bidiagonal matrix whose nonzero entries are \(T_{j,j} = f(q_{n,j})\), \(T_{j,j-1} = -f(q_{n,j-1})\),
$\tilde{H}_n = \textup{diag}( \frac{1}{[\Delta p]_1} - \frac{1}{2n[\Delta p]_1^2},\ldots,  \frac{1}{[\Delta p]_{k}} - \frac{1}{2n[\Delta p]_{k}^2} ) + (\frac{1}{1-[\Delta p]_{k+1}}-\frac{1}{2n(1-[\Delta p]_{k+1})^2}) \bm{1}\bm{1}^T$, $\tilde{b}_n =   ( \frac{1}{2n} - \frac{\sqrt{[\Delta p]_{1}}}{2n^{3/2} [\Delta p]_{k+1}} , \ldots, \frac{1}{2n} - \frac{\sqrt{[\Delta p]_{k}}}{2n^{3/2} [\Delta p]_{k+1}})^T $, and $C_{n,1}$ is a real value that does not depend on $\theta$.
By the definition of $ l(n\Delta F_n(\theta),\Delta \alpha(\theta)) $, we have
\bea 
\log f_n(\theta) &=&  - \frac{n}{2}(\theta - q_n)^T T_n^T\tilde{H}_n T_n(\theta - q_n)  +C_{n,2} +  nR_{n,1}(\theta) +  R_{n,2}(\theta)  ,\\
R_{n,2}(\theta)   &= &  \log \Big[\prod_{j=1}^{k+1}\frac{\Gamma\left( [\Delta \alpha(\theta)]_j \right)  \Gamma( n[\Delta F_n(\theta  )]_j+1)}{ \Gamma\left( [\Delta\alpha(\theta)]_j + n[\Delta F_n(\theta )]_j \right) }\Big],
\eea 
where $C_{n,2}$ is a real value that does not depend on $\theta$.

We have 
\bea 
|nR_{n,1}(q_n+x/\sqrt{n}) -nR_{n,1}(q_n)| &\lesssim&
||x||^3/\sqrt{n} \\
&&+ ||x||(\frac{1}{\sqrt{n}}  +  || \Delta G(q_n+x/\sqrt{n},\theta_0)|| + || \Delta G(q_n,\theta_0)||  ) \\
&&+ || \Delta G(q_n+x/\sqrt{n},\theta_0)||^2 + || \Delta G(q_n,\theta_0)||^2.
\eea        
By Lemma \ref{lemma:AsyEquicontProof}, $|nR_{n,1}(q_n+x/\sqrt{n}) -nR_{n,1}(q_n)|\to 0 $ for fixed $x$ and 
\bea 
|nR_{n,1}(q_n+x/\sqrt{n}) -nR_{n,1}(q_n)| \lesssim \frac{||x||^3}{\sqrt{n}} + x^Tc_n + d_n ,\forall x \in  B_{\sqrt{n}\epsilon} (0 ),
\eea 
for some $c_n$ and $d_n$ with  $||c_n|| \stackrel{p}{\to}  0$ and $d_n \stackrel{p}{\to}  0$.

We show $R_{n,2}(q_n+x/\sqrt{n}) -R_{n,2}(q_n):= \sum_{j=1}^{k+1}T_j \to 0 $ for fixed $x$, where \bean 
T_j\nonumber
&=& \log  \Gamma\left( [\Delta \alpha(q_n+x/\sqrt{n})]_j \right) - \log \Gamma\left( [\Delta \alpha(q_n)]_j \right)  \label{eq:Tj1}\\
&&+ \log ( n[\Delta F_n(q_n+x/\sqrt{n})]_j +1 ) - \log ( n[\Delta F_n(q_n)]_j +1 )\label{eq:Tj2}\\
&&+ ([\Delta \alpha(q_n+x/\sqrt{n})]_j -\frac{1}{2}) \log ( [\Delta \alpha(q_n+x/\sqrt{n})]_j +n [\Delta F_n(q_n+x/\sqrt{n})]_j ) \label{eq:Tj3}\\
&& - ([\Delta \alpha(q_n)]_j -\frac{1}{2}) \log ( [\Delta \alpha(q_n)]_j +n [\Delta F_n(q_n)]_j )  \label{eq:Tj4}\\
&&+ n[\Delta F_n(q_n+x/\sqrt{n})]_j \log \Big( \frac{n[\Delta F_n(q_n+x/\sqrt{n})]_j + 1}{n[\Delta F_n(q_n+x/\sqrt{n})]_j + [\Delta\alpha (q_n+x/\sqrt{n})]_j} \Big) \label{eq:Tj5} \\
&&- n[\Delta F_n(q_n)]_j \log \Big( \frac{n[\Delta F_n(q_n)]_j + 1}{n[\Delta F_n(q_n)]_j + [\Delta\alpha (q_n)]_j} \Big)  \label{eq:Tj6}\\
&&+ R_{n,3}/12\nonumber ,
\eean 
where 
\bea 
\left| R_{n,3} \right| &\leq & \frac{1}{  [n\Delta F_n(q_n)]_j + 1} +\frac{1}{ [n\Delta F_n(q_n)]_j +  [\Delta\alpha (q_n)]_j} \\
&&\frac{1}{  [n\Delta F_n(q_n+x/\sqrt{n})]_j + 1} +\frac{1}{ [n\Delta F_n(q_n+x/\sqrt{n})]_j +  [\Delta\alpha (q_n+x/\sqrt{n})]_j},
\eea
which is obtained by Lemma \ref{lem:logl}.

Since $\Gamma\left( [\Delta \alpha(q_n+x/\sqrt{n})]_j \right)\to \Gamma([\Delta\alpha(\theta_0)]_j)$, $\Gamma\left( [\Delta \alpha(q_n)]_j \right)\to \Gamma([\Delta\alpha(\theta_0)]_j)$ and $[\Delta F_n(q_n)]_j \stackrel{p}{\to} [\Delta F(\theta_0)]_j\neq 0$, we have
\bea 
\eqref{eq:Tj1} &\to& 0,\\
\eqref{eq:Tj2} &=& \log ( 1 + \frac{[\Delta F_n(q_n+x/\sqrt{n})]_j -[\Delta F_n(q_n)]_j }{ [\Delta F_n(q_n)]_j +1/n }) \\
&\stackrel{p}{\to} & 0~ \textup{as }n\to \infty .
\eea 
We have
\bean 
&&\eqref{eq:Tj3}+\eqref{eq:Tj4} \nonumber\\
&=& ([\Delta \alpha(q_n+x/\sqrt{n})]_j -[\Delta \alpha(q_n)]_j)\nonumber\\
&&\times \log ( [\Delta \alpha(q_n+x/\sqrt{n})]_j +n [\Delta F_n(q_n+x/\sqrt{n})]_j ) \label{eq:Tj3_1}\\
&&+ ([\Delta \alpha(q_n)]_j )\nonumber \\
&&\times \log \Big[1 + \frac{([\Delta \alpha(q_n+\frac{x}{\sqrt{n}})]_j -[\Delta \alpha(q_n)]_j)/n + \Delta F_n(q_n+\frac{x}{\sqrt{n}})_j  - [\Delta F_n(q_n)]_j}{[\Delta \alpha(q_n)]_j/n + [\Delta F_n(q_n)]_j} \Big] \label{eq:Tj3_2}.
\eean 
We have $\eqref{eq:Tj3_2}\to 0$ because $[\Delta F_n(q_n)]_j \stackrel{p}{\to} [\Delta F(\theta_0)]_j\neq 0$.
By the Lipschitz continuity assumption of $\alpha$,
\bean\label{eq:ineqTj3_1} 
| \eqref{eq:Tj3_1}   | \lesssim \frac{||x||_2}{\sqrt{n}}  \log (A+n) \to 0.
\eean

Since 
\bea 
\Big( \frac{n\Delta F_n(q_n+x/\sqrt{n})_j + 1}{n\Delta F_n(q_n+x/\sqrt{n})_j + \Delta\alpha (q_n+x/\sqrt{n})_j} \Big)^{n\Delta F_n(q_n+x/\sqrt{n})_j}  \to \exp(1) ~\textup{as }n\to\infty, 
\eea 
we have $\eqref{eq:Tj5} \to 1$. Likewise, $\eqref{eq:Tj6} \to -1$.
Therefore, when $x$ is fixed
\bea 
R_{n,2}(q_n+x/\sqrt{n}) -R_{n,2}(q_n)  \to 0, \textup{as }n\to \infty.
\eea

Next, we show 
\bean\label{Rn2uniform} 
R_{n,2}(q_n+x/\sqrt{n}) -R_{n,2}(q_n) \le \frac{||x||^3}{\sqrt{n}}+  d_n , ~x\in B_{\sqrt{n}\epsilon}(0),
\eean 
for all sufficiently large $n$, where $\limsup_{n\to\infty} d_n <\infty$.

Since $  \inf_{x\in B_{2\sqrt{n}\epsilon}(\theta_0) }\Delta \alpha(\theta+x/\sqrt{n})$ is bounded away from $0$ by the assumption, we have \bea \limsup_{n\to\infty}\sup_{x\in B_{\sqrt{n}\epsilon}(0)}\eqref{eq:Tj1} <\infty.\eea
We have 
\bea 
\limsup_{n\to\infty}\sup_{x\in B_{\sqrt{n}\epsilon}(0)}\eqref{eq:Tj2} &\le& \log \{1+ 2/( [\Delta F_n(q_n)]_j +1/n )\}\to 0,\\
\limsup_{n\to\infty}\sup_{x\in B_{\sqrt{n}\epsilon}(0)}\eqref{eq:Tj3_2}
&\le& \log \{1+ 2/( [\Delta F_n(q_n)]_j + [\Delta \alpha(q_n)]_j )\}\to 0.
\eea 
We have
\bea
\frac{||x||}{\sqrt{n}} \log (A+n) \le \begin{cases}
 \frac{||x||^3}{\sqrt{n}}, &  ~||x|| > \log (A+n)   \\
 \frac{\log (A+n)^2}{\sqrt{n}}&  ~||x|| \le \log (A+n) 
\end{cases}.
\eea 
Thus, by using inquality \eqref{eq:ineqTj3_1}, 
\bea 
\limsup_{n\to\infty}\sup_{x\in B_{\sqrt{n}\epsilon}(0)}\eqref{eq:Tj3_1}  \le  \frac{||x||^3}{\sqrt{n}} + \frac{\log (A+n)^2}{\sqrt{n}}.
\eea 
By the uniform convergence of $F_n$, we have 
\bea 
\limsup_{n\to\infty}\sup_{x\in B_{\sqrt{n}\epsilon}(0)} (\eqref{eq:Tj5} +\eqref{eq:Tj6} ) \to 0.
\eea 
Collecting the results, \eqref{Rn2uniform} is satisfied.

\end{proof}

\subsection{Proof of Proposition \ref{prop:condition2}}\label{ssec:prop2}

For the proof of Proposition \ref{prop:condition2}, we give the following lemma.
\begin{lemma}
\label{lem:negHessianBound}
Let
\[
\log p_{MN}(n x; \, p) 
\;=\;
\log\Gamma(n+1)
\;-\; 
\sum_{j=1}^k \log\Gamma(n\,x_j + 1)
\;+\;
\sum_{j=1}^k n\,x_j \,\log p_j,
\]
where $x = (x_1,\ldots,x_k)\in \mathbb{S}^{k-1}$ and $p=(p_1,\ldots,p_k)\in \mathbb{S}^{k-1}$.  Define 
\bea 
\hat{x}
=
\arg\max_{x \in \mathbb{S}^{k-1}}\,p_{MN}(n x;\,p).
\eea 
Then \(\hat{x}\to p\) as \(n\to\infty\), and in a neighborhood of \(\hat{x}\) 
we have a local quadratic expansion
\[
\log p_{MN}(n x;\,p)
\;=\;
\log p_{MN}\bigl(n \hat{x};\,p\bigr)
-
(\,x_{-k} - \hat{x}_{-k}\,)^T\,H_{n}(x)\,\bigl(x_{-k} - \hat{x}_{-k}\bigr),
\]
where \(H_{n}(x)\) is negative definite.  Moreover, for any $\varepsilon>0$, 
there is a positive constant $C$ (independent of \(n\)) such that
\[
\bigl\lvert H_{n}(x)\bigr\rvert
\;\ge\;
C
\quad
\text{whenever}
\quad
\|\,x - \hat{x}\|\;<\;\varepsilon.
\]
\end{lemma}

\begin{proof}
First, we show that $\hat{x}$ exists uniquely.
Since \(\sum_{j=1}^k -\log\Gamma(n\,x_j+1)\) is strictly concave in the interior 
\(\{x_j>0\}\) and \(\sum_{j} n\,x_j\log p_j\) is linear, $\log p_{MN}(n x;\,p)$ is strictly concave on the interior of the simplex with resepct to $x$. By compactness of \(\mathbb{S}^{k-1}\), a global maximum exists; strict concavity implies uniqueness.
We prove that $\hat{x}$ must lie in the interior (\(\hat{x}_j>0\)) by contradiction.
Suppose that a global maximizer \(\hat{x}=(\hat{x}_1,\dots,\hat{x}_k)\) lies on the boundary of the simplex, say $\hat{x}_j= 0$ for some $j$.
Let $l= \arg\max_i \hat{x}_i$. 
Define $x(\epsilon) = (x_1(\epsilon),\ldots ,x_k(\epsilon))$, where $x_j(\epsilon) = \epsilon$, $x_l(\epsilon) = \hat{x}_l -\epsilon$ and $x_i(\epsilon) =\hat{x}_i$ when $i$ is neither $j$ nor $l$. We have
\bea 
\frac{d}{d\epsilon} \log p_{MN}(n x(\epsilon); p)\big|_{\epsilon=0} = 
-n \psi (1 ) + n(n\hat{x}_l +1) + n\log \frac{p_j}{p_l},
\eea 
where \(\psi(z)=\tfrac{d}{dz}\log\Gamma(z)\) is the digamma function.
Since $\frac{d}{d\epsilon} \log p_{MN}(n x(\epsilon); p)\big|_{\epsilon=0}$ is positive for all sufficiently large $n$. Thus, $\hat{x}$ cannot be the maximizer. 
Furthermore, we show that $\hat{x}$ is bounded away from the boundary: there does not exist a subsequence $a(n)$ and an index $j$ such that $\hat{x}_{a(n),j} \to 0$.
Suppose there exists $j$ and $a(n)$ such that $\hat{x}_{a(n),j} \to 0$.
Since $\hat{x}_{a(n),j}$ is the maximizer on the interior, the following derivative has to be zero:
\bea 
\frac{d}{d\epsilon} \log p_{MN}(n x(\epsilon); p)\big|_{\epsilon=0} &=&
-n \psi (n\hat{x}_{a(n),j} +1 ) + n\psi (n\hat{x}_{a(n),l}  +1) + n\log \frac{p_j}{p_l} .
\eea 
However, since $n\psi (n\hat{x}_{a(n),l}  +1)$ diverges and dominates the other terms, the contraction occurs.

Next, we show that the interior maximizer \(\hat{x}\) of 
\(\log p_{MN}\bigl(n x;\,p\bigr)\) converges to $p$.
Let \(x_{-k} = (x_1,\dots,x_{k-1})\).  Define the function
\[
\ell(x_{-k})
\;=\;
\log p_{MN}\bigl(n x;\,p\bigr)
\;=\;
\log\Gamma(n+1)
\;-\;
\sum_{j=1}^k \log\Gamma\bigl(n\,x_j +1\bigr)
\;+\;
\sum_{j=1}^k n\,x_j \,\log p_j,
\]
where \(x_k=1-\sum_{j=1}^{k-1}x_j\).
Since \(\widehat{x}_{-k}\) is the interior maximum, we have 
\(\nabla \ell(\widehat{x}_{-k})=0\).  In particular, for each \(j=1,\dots,k-1\),
\[
0
\;=\;
\frac{\partial}{\partial x_j}\,\ell(\widehat{x}_{-k})
=
-\,n\,\psi\bigl(n\,\hat{x}_j+1\bigr)
\;+\;
n\,\log p_j
\;-\;
n\,\psi\bigl(n\,\hat{x}_k+1\bigr)\,
\frac{\partial x_k}{\partial x_j}
\;+\;
n\,\log p_k\,\frac{\partial x_k}{\partial x_j},
\]
where \(\psi(z)=\tfrac{d}{dz}\log\Gamma(z)\) is the digamma function, and 
\(\partial x_k/\partial x_j=-1\).  Rearranging gives (for \(j=1,\dots,k-1\)):
\[
-\,\psi\bigl(n\,\hat{x}_j+1\bigr)
\;+\;
\log p_j
\;+\;
\psi\bigl(n\,\hat{x}_k+1\bigr)
\;-\;
\log p_k
\;=\;0.
\]
Thus,
\[
\psi\bigl(n\,\hat{x}_j+1\bigr)
\;=\;
\psi\bigl(n\,\hat{x}_k+1\bigr)
\;+\;
\log\Bigl(\tfrac{p_j}{p_k}\Bigr),
\quad
j=1,\dots,k-1.
\]

We use the asymptotic expansion of the digamma function:
\[
\psi(u)
\;=\;
\log u
\;-\;\frac{1}{2\,u}
\;+\;
O\!\Bigl(\tfrac{1}{u^2}\Bigr).
\]
Since $\hat{x}$ is bounded away from the boundary, we obtain 
\[
\psi\bigl(n\,\hat{x}_j+1\bigr)=
\log\bigl(n\,\hat{x}_j\bigr) + O(1/n),
\quad
\psi\bigl(n\,\hat{x}_k+1\bigr)
=
\log\bigl(n\,\hat{x}_k\bigr)+ O(1/n).
\]
Thus the above equation implies
\bea 
\hat{x}_j  = \hat{x}_k\frac{p_j}{p_k}\exp(O(1/n)) .
\eea 

Since \(\sum_{j=1}^k \hat{x}_j=1\), we get
\[
1
\;=\;
\hat{x}_1 + \cdots + \hat{x}_k
=
\hat{x}_k
\times
\frac{1}{p_k}
\;\bigl(p_1 + \cdots + p_k\bigr) \exp(O(1/n))
\;=\;
\frac{\hat{x}_k}{\,p_k\,}\exp(O(1/n)).
\]
Thus \(\hat{x}_k \to p_k\), which also implies \(\hat{x}_j \to p_j\) 
for all \(j\).

By the 2nd-order Taylor expansion, we have
\[
\log p_{MN}\bigl(n x;\,p\bigr)
\;=\;
\log p_{MN}\bigl(n\hat{x};\,p\bigr)
+\frac{1}{2}
\bigl(x_{-k}-\hat{x}_{-k}\bigr)^T\,H_n(\xi)\,\bigl(x_{-k}-\hat{x}_{-k}\bigr),
\]
where \((k-1) \times (k-1)\) Hessian matrix \( H(x_{-k}) \) has entries:
\[
H_{ij}(x_{-k}) =
\begin{cases}
- n^2 \left[ \psi'(n x_i + 1) + \psi'(n x_k + 1) \right], & i = j, \\
- n^2 \psi'(n x_k + 1), & i \ne j.
\end{cases}
\]
 Because each \(\hat{x}_j>0\) is bounded away from zero, \(\psi'(n\,\hat{x}_j+1)\) 
is bounded below by a constant times \(1/n\), for large \(n\). 
Thus, $\lambda_{\min}(-H_n (x^*))/n \ge C \quad
\text{for all }x\text{ with }\|x-\hat{x}\|<\varepsilon,
$
for some constant \(C>0\).  
\end{proof}

\begin{proof}[Proof of Proposition \ref{prop:condition2}]
Let $\hat{x}_n = argmax_{x \in \mathbb{S}^{k-1} }  \log p_{MN}( n x ; \Delta p)$. By the definition of $f_n$, we have 
    \bean 
    &&\inf_{\theta: |\theta-q_n| >\epsilon} \frac{f_n(q_n) - f_n(\theta)}{n } \nonumber\\
    &\ge &   \inf_{\theta: |\theta-q_n| >\epsilon} \frac{\log p_{MN} ( n \hat{x}_n ;  \Delta p)- \log p_{MN} ( n\Delta F_n(\theta) ;  \Delta p) }{n}  \label{eq:df1}\\
    &&-\frac{| \log p_{MN} ( n \Delta F_n(q_n) ;  \Delta p)- \log p_{MN} ( n \hat{x}_n ;  \Delta p)|}{n}   \label{eq:df2}\\
    &&+\sum_{j=1}^{m+1}\frac{ 
    \log \Gamma([\Delta \alpha(q_n)]_j) + \log \Gamma(n[\Delta F_n(q_n)]_j) -  \log \Gamma([\Delta \alpha(q_n)]_j + n[\Delta F_n(q_n)]_j)
    }{n} ,\label{eq:df4}\\
    &&- \sup_{\theta: |\theta-q_n| >\epsilon}  \sum_{j=1}^{m+1}\frac{ 
    \log \Gamma([\Delta \alpha(\theta)]_j) + \log \Gamma(n[\Delta F_n(\theta)]_j) -  \log \Gamma([\Delta \alpha(\theta)]_j + n[\Delta F_n(\theta)]_j)
    }{n} .\label{eq:df3}
    \eean 
    Since $\hat{x}_n \to \Delta p_{-k}$ (Lemma \ref{lem:negHessianBound}) and $\Delta F_n(q_n) \stackrel{p}{\to } \Delta p_{-k}$, we obtain $\eqref{eq:df2} \stackrel{p}{\to } 0$ and $\eqref{eq:df4} \stackrel{p}{\to } 0$. 

    To show $P(\eqref{eq:df3} \ge  \delta) \to 1$ for arbitrary $\delta>0$,
    we give an upper bound of 
    \bea 
    B_j:=\frac{ 
    \log \Gamma([\Delta \alpha(\theta)]_j) + \log \Gamma(n[\Delta F_n(\theta)]_j) -  \log \Gamma([\Delta \alpha(\theta)]_j + n[\Delta F_n(\theta)]_j)
    }{n},
    \eea considering two cases: (1) $[\Delta F_n(\theta)]_j < M/\sqrt{n}$ and (2) $[\Delta F_n(\theta)]_j \ge M/\sqrt{n}$.
    When $[\Delta F_n(\theta)]_j <  M/\sqrt{n}$, we have 
    \bea 
    B_j&\le& \frac{\log \Gamma(n[\Delta F_n(\theta)]_j)}{n} \\
    &\le& \frac{\log \Gamma( M/\sqrt{n})}{n}\\
    &\le& \frac{M/\sqrt{n} \{\log (M/\sqrt{n}) - 1\} + \log (2\pi)/2 +\sqrt{n}/(12 M  )}{n},
    \eea 
    where the first inequality is satisfied since $\Gamma([\Delta \alpha(\theta)]_j + n[\Delta F_n(\theta)]_j) \ge  \Gamma([\Delta \alpha(\theta)]_j)$ and the third inequality is satisfied by Stirling approximation \eqref{eq:Stirling}.
    Note that the upper bound does not depend on $\theta$ and converges to $0$.
    
    When $[\Delta F_n(\theta)]_j \ge  M/\sqrt{n}$, we have  
    \bea 
    B_j &\le& 
    \frac{n[\Delta F_n(\theta)]_j \log \Big(\frac{n[\Delta F_n(\theta)]_j+1}{n[\Delta F_n(\theta)]_j+ [\Delta \alpha(\theta)]_j} \Big) +\frac{1}{2} \log ( n[\Delta F_n(\theta)]_j +1)}{n} \\
    &&- \frac{([\Delta \alpha(\theta)]_j-\frac{1}{2} ) \log (n[\Delta F_n(\theta)]_j + [\Delta \alpha(\theta)]_j) }{n} \\
    &&+ \frac{([\Delta \alpha(\theta)]_j-\frac{1}{2} ) \log ( [\Delta \alpha(\theta)]_j) - [\Delta \alpha(\theta)]_j +\frac{1}{2} \log (2\pi)}{n} \\
    &&+ \frac{1}{12n} \Big( \frac{1}{ [\Delta \alpha(\theta)]_j} + \frac{1}{ [\Delta \alpha(\theta)]_j +n[\Delta F_n(\theta)]_j  } +\frac{1}{ n[\Delta F_n(\theta)]_j+1 }\Big),\\
    &\le&  \frac{n[\Delta F_n(\theta)]_j}{n}\log \Big(1 +  \frac{1}{n[\Delta F_n(\theta)]_j} \Big)  + \frac{\log (n+1)}{2n} \\
    &&+ \frac{A \{\log (n+A) + \log (A)\}}{n} +\frac{\log (2\pi)}{2n} +\frac{1}{12n}\Big(1 + \frac{2\sqrt{n}}{M'} \Big)\\
    &\le&  \frac{M_2}{n} + \frac{\log (n+1)}{2n} + \frac{A \{\log (n+A) + \log (A)\}}{n} +\frac{\log (2\pi)}{2n} +\frac{1}{12n}\Big(1 + \frac{2\sqrt{n}}{M'} \Big),
    \eea 
    for some positive constant $M_2$, 
    where the second inequality is satisfied since $ [\Delta \alpha(\theta)]_j >M'/\sqrt{n}$ by assumption and 
    $\log \Big(1 +  \frac{1- [\Delta \alpha(\theta)]_j}{n[\Delta F_n(\theta)]_j+ [\Delta \alpha(\theta)]_j} \Big)  \le \log \Big(1 +  \frac{1}{n[\Delta F_n(\theta)]_j} \Big) $, and the last inequality is satisfied because the function $x\log (1+1/x)$ is bounded.
    The upper bound does not depend on $\theta$ and converges to $0$.
    Thus, collecting the first and second cases, we obtain $P(\eqref{eq:df3} \ge  \delta)\to 1$.

    Finally, we show  
    \bea 
    P[ \eqref{eq:df1}  >C] \to 1
    \eea 
    for some positive constant $C$.
    Since  $|\theta-q_n| > \epsilon \Longrightarrow |\theta_{j} - q_{n,j}| >\epsilon/m ~\textup{for some } j$,
    we have 
    \bea 
    |F_n(\theta_j) - \hat{x}_{n,j} | &\ge&  
    |F(\theta_j ) - F(q_{n,j}) | - |F_n(\theta_j) - F(\theta_j ) | - 
    | F(q_n) -   \hat{x}_{n,j} |  \\
     &\ge  & f(q_{n,j}^*) \epsilon/m  - (|F_n(\theta_j) - F(\theta_j ) | + 
    | F(q_{n,j}) -   \hat{x}_{n,j} |),
    \eea 
    where $q_{n,j}^*$ is between $\theta_j$ and $q_{n,j}$.
    Since $|F(q_{n,j}) -   \hat{x}_{n,j}| \to 0$, $|F_n(\theta_j) - F(\theta_j ) | \stackrel{p}{\to } 0$ and $f(q_{n,j}^*) \to f(\theta_0)$, 
    $P(   |F_n(\theta_j) - \hat{x}_{n,j} |  > \epsilon/(2m)) \to 1 $.
    Thus, 
    \bea 
    P\Big[ \eqref{eq:df1}  &\ge& 
    \inf_{\theta: | \Delta F_n(\theta) - \hat{x}_n | >\epsilon/(2m)} \frac{\log p_{MN} ( n \hat{x}_n ;  \Delta p)- \log p_{MN} ( n\Delta F_n(\theta) ;  \Delta p) }{n} \\
    &\ge& \inf_{x: | x - \hat{x}_n | >\epsilon/(2m)} \frac{\log p_{MN} ( n \hat{x}_n ;  \Delta p)- \log p_{MN} ( n x ;  \Delta p) }{n}\Big] \to 1    \eea 
    By Lemma \ref{lem:negHessianBound}, 
    \bea 
\inf_{x: | x - \hat{x}_n | >\epsilon^*} \frac{\log p_{MN} ( n \hat{x}_n ;  \Delta p)- \log p_{MN} ( n x ;  \Delta p) }{n} \ge C
    \eea 
    for some positive constant, which gives $P(\eqref{eq:df1}   >C ) \to 1$.
    By collecting the convergence results of \eqref{eq:df2}-\eqref{eq:df3}, the proof is completed.  
\end{proof}

\bibliographystyle{dcu} 
\bibliography{cdp}